\documentclass[reqno,10pt]{amsart}
\usepackage{amscd,amssymb}
\usepackage{latexsym}
\usepackage[all]{xy}

\def\sreg{{\rm sreg}}
\def\ssing{{\rm ssing}}

\def\into{\hookrightarrow}

\def\Ra{\Rightarrow}

\def\toisom{\widetilde{\to}}
\def\otisom{\widetilde{\leftarrow}}

\def\.{,\dots ,}
\def\wt{\widetilde}
\def\wh{\widehat}
\def\ol{\overline}

\def\Supp{{\rm Supp}}

\def\Spf{{\rm Spf}}
\def\Sp{{\rm Sp}}

\def\Ann{{\rm Ann}}

\def\Spec{{\rm Spec}}
\def\Proj{{\rm Proj}}
\def\Frac{{\rm Frac}}

\def\Bl{{\rm Bl}}

\def\hatBl{{\rm \wh{Bl}}}
\def\red{{\rm red}}

\def\dim{{\rm dim}}

\def\reg{{\rm reg}}
\def\rig{{\rm rig}}
\def\ad{{\rm ad}}
\def\an{{\rm an}}

\def\sing{{\rm sing}}

\def\Id{{\rm Id}}

\def\inv{{\rm inv}}

\def\cha{{\rm char}}

\def\trdeg{{\rm tr.deg.}}

\def\ord{{\rm ord}}

\def\bfA{{\bf A}}

\def\bfF{{\bf F}}

\def\bfN{{\bf N}}

\def\bfP{{\bf P}}
\def\bfQ{{\bf Q}}

\def\bfZ{{\bf Z}}

\def\gtI{{\mathfrak I}}
\def\gtJ{{\mathfrak J}}

\def\gtL{{\mathfrak L}}

\def\gtP{{\mathfrak P}}

\def\gtS{{\mathfrak S}}
\def\gtT{{\mathfrak T}}
\def\gtU{{\mathfrak U}}

\def\gtW{{\mathfrak W}}
\def\gtX{{\mathfrak X}}
\def\gtY{{\mathfrak Y}}
\def\gtZ{{\mathfrak Z}}

\def\calA{{\mathcal A}}
\def\calB{{\mathcal B}}
\def\calC{{\mathcal C}}
\def\calD{{\mathcal D}}

\def\calF{{\mathcal F}}

\def\calI{{\mathcal I}}
\def\calJ{{\mathcal J}}
\def\calK{{\mathcal K}}
\def\calL{{\mathcal L}}
\def\calM{{\mathcal M}}

\def\calO{{\mathcal O}}
\def\calP{{\mathcal P}}

\def\calT{{\mathcal T}}
\def\calU{{\mathcal U}}

\def\calX{{\mathcal X}}
\def\calY{{\mathcal Y}}
\def\calZ{{\mathcal Z}}

\def\oX{{\ol X}}

\def\oZ{{\ol Z}}

\def\of{{\ol f}}

\def\tilS{{\wt S}}
\def\tilT{{\wt T}}

\def\tilW{{\wt W}}

\def\tilY{{\wt Y}}
\def\tilZ{{\wt Z}}

\def\tilg{{\wt g}}

\def\hatA{{\wh A}}
\def\hatB{{\wh B}}

\def\hatX{{\wh X}}

\def\hatZ{{\wh Z}}

\def\hatf{{\wh f}}

\def\ogtJ{{\ol\gtJ}}

\def\ogtT{{\ol\gtT}}

\def\ogtX{{\ol\gtX}}

\def\Kcirc{K^\circ}
\def\Lcirc{L^\circ}

\def\hatOmega{{\wh\Omega}}

\def\alp{{\alpha}}

\def\R+*{{\bf R^*_+}}

\newtheorem{theorsect}{Theorem}[section]
\newtheorem{propsect}[theorsect]{Proposition}

\newtheorem{theor}{Theorem}[subsection]
\newtheorem{prop}[theor]{Proposition}
\newtheorem{lem}[theor]{Lemma}
\newtheorem{cor}[theor]{Corollary}

\newtheorem{question}[theor]{Question}

\theoremstyle{definition}
\newtheorem{definsect}[theorsect]{Definition}

\newtheorem{defin}[theor]{Definition}
\newtheorem{rem}[theor]{Remark}

\begin{document}

\title[Desingularization of quasi-excellent $\bfQ$-schemes]{Desingularization of quasi-excellent schemes
in characteristic zero}
\author{Michael Temkin}
\thanks{The author wishes to thank V. Berkovich and I. Tyomkin for useful discussions and suggestions, E. Bierstone
and P. Milman for communicating the proof of theorem \ref{BMth}, and the referee for careful reading the
manuscript and making many suggestions which improved the exposition. Some inspiration was drawn by the author
from Raynaud's theory which exploits formal blow ups to establish a bridge between the theories of formal
schemes and rigid spaces. A large part of this work was done during author's stay at Max Planck Institute for
Mathematics whose hospitality is gratefully appreciated by the author. }
\address{\tiny{Department of Mathematics, University of Pennsylvania, Philadelphia, PA 19104, USA}}
\email{\scriptsize{temkin@math.upenn.edu}}

\begin{abstract}
Grothendieck proved in EGA IV that if any integral scheme of finite type over a locally noetherian scheme X
admits a desingularization, then X is quasi-excellent, and conjectured that the converse is probably true. We
prove this conjecture for noetherian schemes of characteristic zero. Namely, starting with the resolution of
singularities for algebraic varieties of characteristic zero, we prove the resolution of singularities for
noetherian quasi-excellent Q-schemes.
\end{abstract}

\keywords{Resolution, singularities, quasi-excellent, desingularization}

\maketitle

\section{Introduction}

For a noetherian scheme $X$, let $X_{\rm reg}$ denote the regular locus of $X$. The scheme $X$ is said to {\it
admit a resolution of singularities} if there exists a blow-up $X'\to X$ with center disjoint from $X_{\rm reg}$
and regular $X'$. More generally, for a closed subscheme $Z\subset X$, let $(X,Z)_{\rm reg}$ denote the set of
points $x\in X_{\rm reg}$ such that etale-locally $Z$ is the zero locus of an element $t_1^{n_1} \cdot \dots
\cdot t_d^{n_d}$, where $t_1,\dots,t_d$ is a regular system of parameters. (For example, $(X,\emptyset)_{\rm
reg} = X_{\rm reg}$, and $(X,Z)_\reg=X$ for any regular $X$ with a normal crossing divisor $Z$.) A {\it strict
desingularization (resp. a desingularization) of a pair} $(X,Z)$ is a blow-up $f:X'\to X$ with center disjoint
from $(X,Z)_{\rm reg}$ (resp. from $X_\reg\cap Z_\reg$ and $X_\reg\setminus Z$) and $(X',Z')_{\rm reg}= X'$,
where $Z'=Z\times_X X'$. If, in addition, $f$ is a succession of blow-ups with regular centers, it is said to be
a {\it successive desingularization}. The scheme $X$ is said to admit an {\it embedded} (resp. {\it successive
embedded}) {\it resolution of singularities} if, for any closed subscheme $Z\subset X$, the pair $(X,Z)$ admits
a desingularization (resp. successive desingularization). We remark that usually one does not study strict
desingularizations, but it seems to be a natural extra-condition.

In his celebrated paper \cite{Hir} published in 1964, Hironaka proved that any integral scheme of finite type
over a local quasi-excellent ring of residue characteristic zero admits a successive embedded resolution of
singularities. Recall that a noetherian ring $A$ is said to be quasi-excellent if for any prime ideal
$\wp\subset A$ the canonical homomorphism $A_\wp \to \widehat A_\wp$ is regular, and for any finitely generated
$A$-algebra $B$ ${\rm Spec}(B)_{\rm reg}$ is open in ${\rm Spec}(B)$. (Excellent rings are those which, in
addition to the above two properties, are universally catenary.)

The result of Hironaka is extremely important and has many applications, but its proof is very difficult and
long. It is therefore very natural that mathematicians are still trying to understand and simplify the proof.
Simplified proofs of successive embedded resolution of singularities for integral schemes of finite type over a
field of characteristic zero were first given by Villamayor in \cite{Vil} and Bierstone-Milman in \cite{BM}. In
addition, their desingularization is functorial with respect to smooth morphisms (and strictness of the
desingularization can also be obtained via an additional argument communicated to the author by Bierstone and
Milman, see theorem \ref{BMth}). The works of Villamayor and Bierstone-Milman in their turn served as a basis
for a new generation of proofs, see, for example, \cite{Wl}, \cite{Kol}, \cite{EH}.

On the other hand, Grothendieck proved in \cite[7.9.5]{EGAIV} that if $X$ is a locally noetherian scheme such
that every integral scheme of finite type over $X$ admits a resolution of singularities, then $X$ is
quasi-excellent (i.e. it has a covering by open affine subschemes which are spectra of quasi-excellent rings).
Furthermore, in \cite[7.9.6]{EGAIV} he conjectured that the converse implication is also true, and claimed that
Hironaka's proof gives also resolution of singularities for arbitrary quasi-excellent schemes with residue
fields of characteristic zero, but, as far as we know, Grothendieck's claim was never checked in published
literature.

The purpose of the paper is to show that the existence of embedded resolution of singularities over any
quasi-excellent scheme with residue fields of characteristic zero follows from the corresponding fact for
integral schemes of finite type over fields of characteristic zero. Together with the papers cited above, this
gives a simplified proof of resolution of singularities for arbitrary quasi-excellent schemes with residue
fields of characteristic zero.

In comparison with Hironaka's results, we do not treat successive embedded desingularization, although hope that
this can be done using methods of this paper. On the other hand, we show that Hironaka's theorem for integral
schemes of finite type over a local quasi-excellent ring implies the result stated in Grothendieck's claim
rather easily, see Proposition \ref{locdesprop} and Theorem \ref{Groclaim}.

\subsection{Overview of the paper}

The following theorem is the main result of this paper. We will deduce it from \cite{BM}, and explain in the
appendix how other works can be used instead. Unfortunately, we cannot prove the general strict desingularization
(the bottleneck being proposition \ref{regprop}, see remark \ref{basarem}), so we weaken it as follows: a
desingularization $f:X'\to X$ of a pair $(X,Z)$ is {\em semi-strict} if the center of $f$ is disjoint from the
subset of $(X,Z)_\reg$ in which the irreducible components of $Z$ have no self-intersections (it is the strictly
monomial locus of $Z$ in $X$ in the sense of definition \ref{monomdef}, and it is denoted $(X,Z)_\sreg$ in the
paper starting with definition \ref{reglocdef}).

\begin{theorsect}
\label{mainth} Let $X$ be a noetherian scheme of characteristic zero, then the following conditions are
equivalent:

(i) $X$ is quasi-excellent;

(ii) any integral scheme of finite type over $X$ admits a desingularization;

(iii) any integral scheme of finite type over $X$ admits a semi-strict embedded resolution of singularities.
\end{theorsect}

Since (iii) is obviously stronger than (ii), and the implication (ii)$\Rightarrow$(i) is due to Grothendieck,
the theorem is equivalent to proving that any integral quasi-excellent scheme admits semi-strict embedded
resolution of singularities. Now, let us give a more detailed description of the paper.

We study basic properties of blow ups and desingularizations in \S\ref{localch}, and the main result is
proposition \ref{locdesprop} which states that there is resolution of singularities over a noetherian
quasi-excellent scheme $k$ if and only if any scheme $Y$ isomorphic to a blow up of a local $k$-scheme of
essentially finite type admits a desingularization. Thus, up to a not so difficult proposition \ref{locdesprop},
Hironaka's theorem implies that any noetherian quasi-excellent scheme admits a desingularization. In particular,
we obtain the equivalence of (i) and (ii) in the main theorem.

The direct implication in \ref{locdesprop} is straightforward, and the opposite one is proved by a simple
argument which is used few more times in the paper. Therefore we outline it here, assuming for simplicity that
$Z$ is empty. Consider a scheme $X$ with a subset $S\subset X$, and let $X'\to X$ be a blow up. Let us say that
$f$ desingularizes $X$ {\em over} $S$ if $f^{-1}(S)\subset X'_\reg$. We start with the identity morphism $\Id_X$
which desingularizes $X$ over $S_0=X_\reg$ and construct a desingularization of $X$ using noetherian induction.
The induction step is as follows: we start with a blow up $X'\to X$ desingularizing $X$ over an open set
$S\subset X$, choose a maximal (or generic) point $x$ of the complement of $S$ and note that by our assumptions,
the scheme $X'_x=\Spec(\calO_{X,x})\times_X X'$ (which is a blow up of $\Spec(\calO_{X,x})$) admits a
desingularization $f'_x:X''_x\to X'_x$. Then we extend $f'_x$ to a blow up $f':X''\to X'$ trivial over
$f^{-1}(S)$ and note that the composition $f\circ f':X''\to X$ desingularizes $X$ over an open set $S'$
containing $S$ and $x$. Note that the extension $f'$ of $f'_x$ can be extremely bad above proper specializations
of $x$; in particular, the resulting desingularization can be not successive even when $f'_x$ is a successive
one.

\S\ref{localch} is organized as follows. In \S\ref{blowsec}, we study extensions of ideals and blow ups, and
introduce formal blow ups. In the next section, we fix our desingularization terminology, and we prove
proposition \ref{locdesprop} in \S\ref{ressec}.

In \S\ref{formch}, we prove the equivalence (i)$\Leftrightarrow$(ii) once again, but this time using only
desingularization of integral schemes of finite type over a field of characteristic zero. The main idea is as
follows: one can construct a desingularization of a quasi-excellent scheme $X$ from the desingularization of its
completion $\gtX$ along the singular locus, and if the latter is algebraizable, i.e. is isomorphic to a
completion of a scheme $\calX$ of finite type over a field, then it suffices to know how to desingularize
$\calX$. For example, one can desingularize isolated singularities because any complete local ring with an
isolated singularity is algebraizable due to Artin, see \cite[3.8]{Art}. The results of Artin were generalized
by Elkik in \cite{Elk}, in particular, she proved that any affine rig-smooth formal scheme of finite type over a
complete ring with a principal ideal of definition is algebraizable. Since rig-smoothness is equivalent to
rig-regularity in the characteristic zero case, Elkik's results can be applied to desingularize an affine
quasi-excellent scheme $X$ whose singular locus $X_\sing$ is of finite type over a field of characteristic zero.
The case of an arbitrary $X$ is obtained from this one by using a noetherian induction argument similarly to the
proof of proposition \ref{locdesprop}.

\S\ref{formch} is structured as follows. We introduce quasi-excellent formal schemes in \S\ref{formregsec} and
define for them notions of regularity, reducedness, etc. In next section we introduce an important class of {\em
special} formal schemes which are quasi-excellent by results of Valabrega. Then, in \S\ref{rigsmoothsec}, we
deduce from Elkik's theorem that certain rig-smooth special formal schemes are algebraizable, see proposition
\ref{algebprop}. It is known that our proposition is a form of Elkik's theorem, but we prefer to prove it
because of lack of reference. Finally, in \S\ref{formmainsec}, we use proposition \ref{algebprop} to
desingularize certain rig-smooth special formal schemes, see theorem \ref{formtheor}. Desingularization of
quasi-excellent schemes follows easily.

We cannot treat embedded desingularization in \S\ref{formch} because Elkik's theorem algebraizes certain formal
schemes, but not pairs consisting of a formal scheme and a divisor. Although the author expects that one can
algebraize certain rig-monomial divisors (thus generalizing Elkik's theorem), this question is not studied in
the paper. We prove proposition \ref{regprop} instead, and use it to monomialize strict transform of a divisor.
Combining this proposition with the results of \S\ref{formch}, we are able to prove theorem \ref{mainth} in
general. At the end of \S\ref{strch} we desingularize formal rig-regular schemes, see theorem \ref{formmainth}.

The paper contains an appendix where we study a connection between semi-strict embedded desingularization and
standard desingularization results in which the entire set $(X,Z)_\sing\cup Z_\sing$ can be modified (for
example, it happens when one applies Main Theorem II of \cite{Hir}). We prove in the appendix that one can
deduce semi-strict embedded desingularization from Main Theorem II of Hironaka or its analogs.

\subsection{Terminology and notation}
If $X$ is a scheme with a closed subscheme $Z$, then $|Z|$ denotes the {\em support} of $Z$, i.e. the underlying
set of $Z$ considered as a closed subset of $X$. The support $\Supp(\calI)$ of an ideal $\calI\subset\calO_X$ is
the support of the associated closed subscheme $Z=\Spec(\calO_X/\calI)$. We will pass freely from reduced closed
subschemes to closed subsets and vice versa, but we will use different notation: we write $Z\subset Z'$, $Z\cap
Z'$, $Z\setminus Z'=Z-(Z\cap Z')$ and $f^{-1}(Z)$ (where $f:X'\to X$ is a morphism) when working with subsets,
and we write $Z\into Z'$, $Z\times_X Z'$ and $Z\times_X X'$ when working with subschemes. By $X^c$, $X^{\le c}$,
etc., we denote the sets of points of codimension $c$, of codimension at most $c$, etc. In particular, $X^0$ is
the set of maximal points of $X$.

Recall that a noetherian ring $A$ is called {\em quasi-excellent} if for any prime ideal $p\subset A$, the
completion morphism $\phi:A_p\to\hatA_p$ is {\em regular} (i.e. $\phi$ is flat and has geometrically regular
fibers) and for any finitely generated $A$-algebra $B$, the regular locus of $\Spec(B)$ is open. A universally
catenary quasi-excellent ring is called {\em excellent}. A scheme $X$ is called {\em (quasi-) excellent} if it
admits an open covering by spectra of (quasi-) excellent rings.

If $A$ is a ring with an ideal $P$ and an $A$-ring $B$, then by {\em $P$-adic completion} $\hatB_P$ of $B$ we
mean the separated completion of $B$ in $(PB)$-adic topology. We say that $B$ is {\em $P$-adic} if
$B\toisom\hatB_P$. Similarly, if $X'\to X$ is a morphism of schemes and $\calI\subset\calO_X$ is an ideal, then
$\hatX'_\calI$ denotes the $\calI$-adic completion of $X'$, i.e. the formal completion of $X'$ along
$\Spec(\calO_{X'}/\calI\calO_{X'})$.

We refer to \cite[\S1.10]{EGAI} for basic properties of formal schemes. Any formal scheme appearing in this
paper is automatically assumed to be locally noetherian. Moreover, if not said to the contrary, it is assumed to
be noetherian. Given a locally noetherian formal scheme $\gtX$, its {\em closed fiber} $\gtX_s$ is defined as
$\Spec(\calO_\gtX/\gtP)$, where $\gtP$ is the biggest ideal of definition. Recall that $\gtX_s$ is a reduced
closed (formal) subscheme of $\gtX$, homeomorphic to it as a topological space. If $\gtX$ is a formal scheme of
finite type over a complete discrete valuation ring, then one can attach to $\gtX$ a generic fiber
$\gtX_\eta^\rig$ (resp. $\gtX_\eta^\an$, resp. $\gtX_\eta^\ad$), which is a rigid (resp. analytic, resp. adic)
space, see \cite[Ch. 9]{BGR} (resp. \cite{Ber1}, resp. \cite{Hub}). Only minimal familiarity with classical
rigid spaces is required in this paper. The author expects, however, that analytic or adic spaces can be useful
in attacking question \ref{algebrquest}. Note also that adic and generalized rigid generic fibers are defined
for arbitrary noetherian formal schemes, see \cite{Hub} or \cite[\S5]{BL1}.

\section{Blow ups and desingularization}
\label{localch} In this section, we establish some properties of blow ups of schemes and formal schemes which
will be used later. The main result is proposition \ref{locdesprop} which localizes the desingularization
problem by reducing desingularization of a general scheme to desingularization of blow ups of local schemes.

\subsection{Ideals and blow ups}
\label{blowsec} Let us consider the following situation which is a particular case of the situation considered
in \cite[8.2.13]{EGAIV}. Assume that $X$ is a noetherian scheme and $S_\alp$ is a filtered family of open
subschemes such that the transition morphisms $S_\beta\to S_\alpha$ are affine. Then there exists an $X$-scheme
$S=\projlim_\alp S_\alp$, the structure morphism $i:S\to X$ maps $S$ homeomorphically onto its image, and
$\calO_S$ is isomorphic to the restriction of $\calO_X$ on $i(S)$. It follows that $S$ is noetherian too. We
will identify $S$ with $(i(S),\calO_X|_{i(S)})$ and say that it is a {\em pro-open pro-subscheme} of $X$. A
typical example is obtained when $S_\alp$ are affine neighborhoods of a point $x$; then
$S\toisom\Spec(\calO_x)$. Another example is obtained from this one by base change with respect to a morphism
$X'\to X$.

\begin{lem}
\label{extendlem} Keep the above notation and assume that we are given an ideal $\calI_S\subset\calO_S$. Let
$Z_S$ denote the support of $\calI_S$ and $Z$ be its Zariski closure in $X$. Then there exists an ideal
$\calI\subset\calO_X$ such that $\calI|_S=\calI_S$ and the support of $\calI$ is $Z$.
\end{lem}
\begin{proof}
Since $S$ is noetherian, there is a one-to-one correspondence between ideals in $\calO_S$, closed subschemes of
$S$ and closed subschemes of finite presentation. For a scheme $Y$, let $\gtS(Y)$ denote the set of closed
subschemes of $Y$. Note that $\injlim_\alp\gtS(S_\alp)\toisom\gtS(S)$ by \cite[8.6.3]{EGAIV}, and the map
$\gtS(X)\to\gtS(S)$ is surjective because the map $\gtS(X)\to\gtS(S_\alp)$ is surjective for any $\alp$ by
\cite[6.9.7]{EGAI}. Hence $\calI_S=\calI'|_S$ for some $\calI'\subset\calO_X$, and the first claim of the lemma
is satisfied.

Let $Z'$ be the support of $\calI'$, then its intersection with $S$ coincides with $Z_S$, and $Z_S$ is closed
under generalizations in $Z'$ because $S$ is closed under generalizations in $X$. It follows that the set $Z'^0$
of maximal points of $Z'$ is a union of $Z_S^0$ and a finite set $Z''^0$ disjoint from $Z_S$. Hence $Z'=Z\cup
Z''$, where $Z''$ is a closed set disjoint from $S$. To finish the proof, we have to "correct" $\calI'$ over
$Z''$. Note that $U:=X\setminus Z''$ is a neighborhood of $S$, and the support of $\calI'|_U$ coincides with
$Z\cap U$. Since $Z\cap U=Z-(Z\cap Z'')$ is closed in $Y:=X-(Z\cap Z'')$, we can trivially extend $\calI'|_U$ to
an ideal $\calI_Y\subset\calO_Y$ on $Y$. Indeed, the sheaves of ideals $\calI'|_U$ and $\calO_{Y\setminus Z}$ on
the open subschemes $U$ and $Y\setminus Z$, respectively, agree over the intersection $U\cap (Y\setminus
Z)=U\setminus Z$, hence they glue to an ideal $\calI_Y$ on $U\cup(Y\setminus Z)=Y$. Finally, let $\calI$ be any
extension of $\calI_Y$ to $X$ (it exists by \cite[6.9.7]{EGAI}), then the support of $\calI$ is contained in
$(Z\cap U)\cup (Z\cap Z'')=Z$, as required.
\end{proof}

Next we discuss sheaves of ideals on a noetherian formal scheme $\gtX$. Let $\gtI\subset\calO_\gtX$ be an ideal
with the associated closed formal subscheme $\gtZ$. It is not clear how to define the reduction of $\gtZ$ in
general, so we are forced to give the following definition. Given two ideals $\gtI,\gtJ$, we say that {\em the
support of $\gtI$ is contained in the support of $\gtJ$} if $\gtJ^n\subset\gtI$ for some $n$. We will not use
the following side remark in our proofs, though it will be mentioned in few more side remarks.

\begin{rem}\label{adicrem}
If one takes into account the generic fiber $\gtX_\eta^\ad$, then one obtains a reasonable set-theoretical
description of support of an ideal on a formal scheme. Indeed, if $\gtI$ is an ideal such that the corresponding
closed subscheme $\gtZ$ is reduced, then $\gtI$ is uniquely determined by its adic support
$|\gtZ^\ad|=|\gtZ|\sqcup|\gtZ_\eta^\ad|\subset|\gtX^\ad|$ (combining the usual support with the "generic" one).
Note that $|\gtZ|$ itself is far too small to determine $\gtI$.
\end{rem}

In general, one cannot extend to $\gtX$ an ideal defined on an open formal subscheme, and one cannot algebraize
an ideal on the formal completion $\gtX=\hatX_Z$ of a scheme $X$ along a closed subscheme $Z$. The situation
with open ideals is better. Any open ideal $\gtI$ defines a closed formal subscheme $\gtZ=\Spf(\calO_\gtX/\gtI)$
which is a scheme (i.e. its ideal of definition is nilpotent) supported on $\gtX_s$. Actually, if $\gtP$ is an
ideal of definition, then $\gtZ$ is a closed subscheme of some $\Spec(\gtX/\gtP^n)$. In particular, we can and
do define the {\em support} of $\gtI$ as a closed subset of $\gtX$. If $\gtX=\hatX_Z$, then the completion
induces a bijective correspondence between ideals $\calI\subset\calO_X$ supported on $|Z|$ and open ideals
$\gtI\subset\calO_\gtX$. In other words, open ideals are algebraizable.

\begin{lem}
\label{formextendlem} Let $\gtX$ be a noetherian formal scheme with an open formal subscheme $\gtY$ and a closed
subset $Z\subset\gtX$. Then any open ideal $\gtI\subset\calO_\gtY$ with support in $Z\cap\gtY$ extends to an
open ideal $\gtJ\subset\calO_\gtX$ with support in $Z$.
\end{lem}
\begin{proof}
Consider the open formal subscheme $\gtY'=\gtY\cup(\gtX\setminus Z)$. The ideal $\gtI$ can be extended trivially
to an open ideal $\gtI'\subset\calO_{\gtY'}$ (as in the proof of lemma \ref{extendlem}, we use that the support
of $\gtI$ is closed in $\gtY'$ to glue $\gtI'$ from $\gtI$ and $\calO_{\gtY'\setminus Z}$). Now it suffices to
find an arbitrary extension of $\gtI'$ to an open ideal $\gtJ\subset\calO_\gtX$. Choose an ideal of definition
$\gtP\subset\calO_\gtX$ such that $\gtP'=\gtP|_{\gtY'}$ is contained in $\gtI'$, and consider the schemes
$X=(\gtX_s,\calO_\gtX/\gtP)$, $Y'=(\gtY'_s,\calO_{\gtY'}/\gtP')$ and the ideal $\calI'=\gtI'\calO_{Y'}$. Then
$\calI'$ extends to an ideal $\calJ\subset\calO_X$ by \cite[6.9.7]{EGAI}, and the preimage
$\gtJ\subset\calO_\gtX$ of $\calJ$ is a required extension of $\gtI'$.
\end{proof}

Next, we recall basic facts about blow ups, see \cite[\S1]{Con1}, for more details. If $X$ is a scheme with a
finitely generated ideal $\calI\subset\calO_X$, then the $X$-scheme
$X'=\Proj(\calO_X\oplus\calI\oplus\calI^2\oplus\dots)$ is $X$-projective and the structure morphism $X'\to X$ is
an isomorphism over $X\setminus\Supp(\calI)$. The pair $(X',\calI)$ is called the {\em blow up of $X$ along
$\calI$} and is denoted $\Bl_\calI(X)$. The ideal $\calI\calO_{X'}$ is invertible and $X'$ is the final object
in the category of $X$-schemes such that the preimage of $\calI$ is invertible. The construction of blow ups
commutes with localizations (and, more generally, with flat base changes). If $X=\Spec(A)$ and
$I=\calI(X)\subset A$, then $X'$ is glued from schemes $X'_g=\Spec(A[\frac{I}{g}])$ with $g\in I$, where
$A[\frac{I}{g}]$ is the subring of $A_g$ generated by $\frac{I}{g}$. Usually the schemes $X'_g$ are called {\em
charts} of the blow up.

As usual, we will omit the ideal in the notation of a blow up, and say simply "a blow up $f:X'\to X$", or even
"a blow up $X'$ of $X$", however, we will take the ideal into account in definitions of $V$-admissibility and
strict transforms. Note that the $X$-scheme $X'=\Bl_\calI(X)$ can be obtained by blowing up other ideals, for
example $\Bl_{\calI^n}(X)\toisom\Bl_\calI(X)$ for any $n>0$. Sometimes we will say that $X'$ is a blow up of $X$
along $Y=\Spec(\calO_X/\calI)$, and $Y$ is the {\em center} of the blow up, and write $X'=\Bl_Y(X)$. For any
open subscheme $V$, the blow up $f:X'\to X$ is called {\em $V$-admissible} if its center is disjoint from $V$.
Sometimes it will be more convenient to express the same property in terms of the complementary closed set, so
we say that $f$ is {\em $T$-supported} for a closed subscheme (or subset) $T\into X$ if $|Y|\subset |T|$ (i.e.
$f$ is $(X\setminus T)$-admissible). More generally, given a morphism $g:X\to S$, a closed subscheme $R\into S$
and an open subscheme $U=S\setminus R$, we say that $f$ is $U$-admissible (or $R$-supported) if $f$ is
$g^{-1}(U)$-admissible. We will make an intensive use of the following well known result. A simple and natural
proof of this fact given in \cite[5.1.4]{RG} is incomplete, and we refer to \cite[1.2]{Con1} for a surprisingly
involved full proof due to Raynaud.

\begin{lem}\label{compblow}
If $X$ is coherent (i.e. quasi-compact and quasi-separated), $V\into X$ is open and $T=X\setminus V$, then a
composition of $V$-admissible (or $T$-supported) blow ups is a $V$-admissible (or $T$-supported) blow up.
\end{lem}

If $V$ is an open subscheme of $X$ such that a blow up $f:X'\to X$ is an isomorphism over $V$, then it still can
happen that $f$ is not isomorphic to a $V$-admissible blow up. For example, it is the case when
$X=\Spec(k[x,y,z,t]/(xy-zt))$, $V=X_\reg$ is the complement of the origin $s$ and $I=(x,y)$ defines a Weil
divisor which is not Cartier. Then $X'=\Bl_I(X)$ is a small resolution of $X$, and it cannot be obtained by
blowing up an ideal supported on $s$ because the preimage of $s$ is not a divisor, but a curve. Nevertheless,
$X'$ is dominated by a $V$-admissible blow up. More generally, we will need the following lemma.

\begin{lem}
\label{admisslem} Let $X$ be a coherent scheme with a schematically dense open subscheme $U$ and $f:X'\to X$ be
a $U$-modification, i.e. a proper morphism such that $f^{-1}(U)$ is schematically dense in $X'$ and is
$X$-isomorphic to $U$. Then there exists a $U$-admissible blow up $X''\to X$ which factors through $X'$.
\end{lem}
\begin{proof}
Apply the flattening theorem of Raynaud and Gruson, see \cite[5.2.2]{RG}, to the morphism $f:X'\to X$ and the
sheaf $\calO_{X'}$ (we set $S=X$, $X=X'$ and $\calM=\calO_{X'}$ in the loc.cit.). By the theorem, there exists a
$U$-admissible blow up $\oX\to X$ such that the following condition holds: let $\of:\oX'\to\oX$ denote the base
change of $f$ and $\calF$ denote the strict transform of $\calO_{X'}$, then $\calF$ is $\calO_\oX$-flat.

Note that $U$ is a schematically dense open subscheme of $X,X'$ and $\oX$, and let $X''$ be the schematic
closure of the image of $U$ under the diagonal morphism $U\to\oX'$. Then $X''$ is the minimal $U$-modification
of $X$ which dominates both $X'$ and $\oX$, and $\calO_{X''}$ is isomorphic to the quotient of $\calO_{\oX'}$ by
the maximal submodule supported on the preimage of $X\setminus U$, i.e. $\calO_{X''}\toisom\calF$. Thus,
$g:X''\to\oX$ is a flat $U$-modification, and we will prove that it is an isomorphism. It will follow then that
$X''$ is a required $U$-admissible blow up which dominates $X'$.

To check that $g$ is an isomorphism we may work locally on $\oX$, so assume that $\oX=\Spec(A)$. Since $g$ is
flat and an isomorphism over $U$, the fibers of $g$ are discrete, i.e. $g$ is quasi-finite. Since $g$ is proper,
it is finite, and therefore $X''=\Spec(B)$. By flatness of $g$, $g_*(\calO_{X''})$ is a locally free
$\calO_\oX$-sheaf. The rank of $g_*(\calO_{X''})$ is $1$ at any point of $\oX$ because it is so on a dense
subscheme $U$. We obtain that $B=hA$ for an element $h\in B$, hence $1=ha$ for some $a\in A$. Moreover, $a$ is
invertible in $A$ because the map $X''\to X$ is surjective, and we obtain that $B=A$ as claimed.
\end{proof}

Let $\gtX$ be a noetherian formal scheme and $\gtI\subset\calO_\gtX$ be an ideal. If $\gtI$ is open, then a
notion of admissible formal blow up along $\gtI$ is defined in \cite[\S2]{BL1}. Our last goal in this section is
to introduce formal blow ups along arbitrary ideals and study their basic properties (in the case of an open
ideal our definition is slightly different because we do not restrict to admissible formal schemes). However,
blow ups along not open ideals will appear in the last section of the paper, and until then the results of
\cite[\S2]{BL1} cover our needs, so the reader can consult loc.cit. instead of reading the rest of this section.

\begin{rem}\label{Nicrem}
Arbitrary formal blow ups (generalizing admissible formal blow ups) were defined independently by J. Nicaise. In
a recent work \cite{Nic} on a trace formula and motivic integration, he establishes some basic properties of
formal blow ups, including our lemma \ref{formblowlem} below (see Proposition 2.16 in loc.cit).
\end{rem}

Assume that $\gtX=\Spf(A)$ is affine, and let $I\subset A$ be the ideal corresponding to $\gtI$ and $P\subset A$
be an ideal of definition. We define the formal blow up $\gtX'=\hatBl_I(A)$ of $\gtX$ along $\gtI$ as the
$P$-adic completion of $X'=\Bl_I(\Spec(A))$. Since $X'$ is glued from affine charts $X'_g=\Spec(A[\frac{I}{g}])$
with $g\in I$, its completion is glued from affine formal schemes $\gtX'_g=\Spf(A\{\frac{I}{g}\})$, where
$A\{\frac{I}{g}\}$ is the $P$-adic completion of $A[\frac{I}{g}]$. Let us give an explicit description of
$A\{\frac{I}{g}\}$. First, we note that the homomorphism $A\{\frac{I}{g}\}\to A_{\{g\}}$ can be not injective
(for example, the target is zero when $P\subset(g)$), so it is of no use for us. From other side, it is well
known that if $I=(f_1\. f_n)$, then $A[\frac{I}{g}]$ can be described as the quotient of the ring $A'=A[T_1\.
T_n]/(gT_1-f_1\. gT_n-f_n)$ by its $g$-torsion. Note that the completion of $A'$ if isomorphic to
$\hatA'=A\{T_1\. T_n\}/(gT_1-f_1\. gT_n-f_n)$. Since $A[\frac{I}{g}]$ is noetherian, its completion
$A\{\frac{I}{g}\}$ is flat over it, and, in particular, $A\{\frac{I}{g}\}$ has no $g$-torsion. It follows that
$A\{\frac{I}{g}\}$ is isomorphic to the quotient of $\hatA'$ by its $g$-torsion.

\begin{lem}
\label{formlem} Let $A$ be a noetherian ring with ideals $I$ and $P$, and $\hatA$ be the $P$-adic completion of
$A$. Then the $P$-adic completion of $\Bl_I(A)$ is canonically isomorphic to $\hatBl_{I\hatA}(\hatA)$.
\end{lem}
\begin{proof}
Let $f_1\. f_n$ be generators of $I$ and $g\in I$ be an element. We have to prove that the $P$-adic completion
of $A[\frac{I}{g}]$ is canonically isomorphic to $\hatA\{\frac{I}{g}\}$. Obviously, the $P$-adic completion of
$A'=A[T_1\. T_n]/(gT_1-f_1\. gT_n-f_n)$ is isomorphic to $\hatA'=\hatA\{T_1\. T_n\}/(gT_1-f_1\. gT_n-f_n)$. By
the same flatness argument as above, the $P$-adic completion of $A'/(g)$-torsion is isomorphic to
$\hatA'/(g)$-torsion.
\end{proof}

If $f\in A$ is an element, then $\hatBl_{IA_{\{f\}}}(\gtX_{\{f\}})$ is isomorphic to the completion of
$\Bl_{IA_f}(A_f)$. Since usual blow ups are compatible with localizations, the latter is isomorphic to the
completion of $(\Bl_I(A))_f$, which in its turn is isomorphic to $(\hatBl_I(\gtX))_{\{f\}}$. We see that formal
blow ups of affine formal schemes are compatible with formal localizations, and it follows, in particular, that
for any locally noetherian formal scheme $\gtX$ with ideal $\gtI\subset\calO_\gtX$, one can define the formal
blow up $\hatBl_\gtI(\gtX)$ by gluing formal blow ups of open affine formal subschemes of $\gtX$. We say that a
formal blow up $X'=\hatBl_\gtI(X)\to X$ is {\em $\gtJ$-supported} if the support of $\gtI$ is contained in the
support of $\gtJ$. Furthermore, if $\gtI$ is open, then its support lies in a closed subset $T\subset\gtX_s$,
and we say then that $X'$ is $T$-supported.

\begin{lem}
\label{formblowlem} Let $X$ be a noetherian scheme with two ideals $\calI,\calP\subset\calO_X$, $\gtX$ be the
$\calP$-adic completion of $X$ and $\gtI=\calI\calO_\gtX$. Then the $\calP$-adic completion of $\Bl_\calI(X)$ is
canonically isomorphic to $\hatBl_\gtI(\gtX)$.
\end{lem}
\begin{proof}
Both formal completion and blow up along a closed subscheme are defined locally on $X$, so it suffice to
consider the affine case, which was established in the previous lemma.
\end{proof}

Let $\gtX$ be a noetherian formal scheme with a closed formal subscheme $\gtT$. If $\gtT$ is a scheme, then
Raynaud proved that a composition of $\gtT$-supported formal blow ups is isomorphic to a $\gtT$-supported formal
blow up, see \cite[2.5]{BL1}. The same is true for general formal blow ups. Let $\gtI\subset\calO_\gtX$ be a
$\gtT$-supported ideal with formal blow up $\hatf:\gtX'=\hatBl_\gtI(\gtX)\to\gtX$ and $\gtJ\subset\calO_{\gtX'}$
be a $\gtT\times_\gtX\gtX'$-supported ideal with formal blow up $\hatf':\gtX''=\hatBl_\gtJ(\gtX')\to\gtX'$.

\begin{lem}
\label{formsupplem} Keep the above notation, then $\gtX''$ is $\gtX$-isomorphic to a $\gtT$-supported blow up.
\end{lem}
\begin{proof}
Set $\gtI'=\gtI\calO_{\gtX'}$. Let $m,n$ be natural numbers, then $(\hatf_*(\gtI'^n\gtJ))^m$ is an ideal in the
$\calO_\gtX$-algebra $\hatf_*(\calO_{\gtX'})$ and its preimage in $\calO_\gtX$ is an ideal $\gtL$ (depending on
$m$ and $n$). We will prove that $\gtX''\to\gtX$ is isomorphic to the blow up of $\gtX$ along $\gtL\gtI$ for
$n>n_0$ and $m>m_0(n)$.

The latter statement can be checked locally on $\gtX$ because $\gtX$ is quasi-compact. So, we can assume that
$\gtX=\Spf(A)$ for a $P$-adic ring $A$. Let $I\subset A$ denote the ideal corresponding to $\gtI$, $X=\Spec(A)$,
$X'=\Bl_I(X)$ and $f:X'\to X$ be the blow up morphism. Since $X'$ is $X$-proper and $\gtX'$ is isomorphic to the
$P$-adic completion of $\gtX$, we can apply Grothendieck's existence theorem, see \cite{EGAIII}, theorem 5.1.4
and corollary 5.1.8, to find an algebraization $\calJ\subset\calO_{X'}$ of $\gtJ$. Set $T=\Spec(A/m)$, where $m$
is such that $\gtT=\Spf(A/m)$, then $I$ is supported on $T$ and $\calJ$ is supported on the preimage of $T$ in
$X'$. By lemma \ref{compblow}, $X''=\Bl_\calJ(X')$ is $X$-isomorphic to a $T$-supported blow up of $X$. Since
$\gtX''$ is isomorphic to the $P$-adic completion of $X''$, we already obtain that $\gtX''$ is isomorphic to a
$\gtT$-supported blow up of $\gtX$. However, as we mentioned above, we have to describe the blow up explicitly,
and this will require a closer look on the proof of \cite[1.2]{Con1}.

Set $\calI'=I\calO_{X'}$. The proof in loc.cit. starts with an observation that for sufficiently large $n$ and
ideal $\calM=\calI'^n\calJ$, the map $f^*f_*(\calM)\to\calM$ is surjective. Then a sufficiently large submodule
$\calK\subset f_*(\calM)$ of finite type is chosen ($X$ can be non-noetherian in loc.cit.). Since $f_*(\calM)$
is coherent in our situation, one can actually choose $\calK=f_*(\calM)$. Finally, one defines
$\calL\subset\calO_X$ in loc.cit. as the preimage of $\calK^m$ under the homomorphism $\calO_X\to
f_*(\calO_{X'})$, and proves that for sufficiently large $m$, $X''\toisom\Bl_{I\calL}(X)$.

By lemma \ref{formblowlem}, $\gtX''\toisom\hatBl_{I\calL}(\gtX)$, so we have only to prove that $\calL=\gtL$ (as
an ideal in $A$). Note that $\hatf_*(\calO_{\gtX'})$ is isomorphic to the $P$-adic completion of
$f_*(\calO_{X'})$ by Grothendieck's theorem on formal functions, see \cite[4.1.5]{EGAIII}, but $f_*(\calO_{X'})$
coincides with its completion because it is a finite $A$-algebra. By the same argument, $\hatf_*(\gtI'^n\gtJ)$
and $f_*(\calI'^n\calJ)$ define the same ideal in $f_*(\calO_{X'})$, and therefore $\gtL=\calL$.
\end{proof}

\subsection{Desingularization of a pair} \label{desingsec}
Let $X$ be a locally noetherian scheme and $Z$ be a closed subscheme. We say that $Z$ is {\em a Cartier divisor}
if it is locally given by a single regular element, i.e. for any point $z\in Z$, there exists a not zero divisor
$f_z\in\calO_{X,z}$ with $\calO_{X,z}/f_z\calO_{X,z}\toisom\calO_{Z,z}$. This condition is equivalent to
requiring that the ideal $\calI\subset\calO_X$ defining $Z$ is invertible.

\begin{defin}\label{monomdef}
We say that $Z$ is a {\em strictly monomial divisor} if $X$ is regular in a neighborhood of $Z$ and for any
point $x\in Z$ there exists a regular sequence of parameters $u_1\. u_d\in\calO_{X,x}$ such that $Z$ is given by
an equation $\prod_{i=1}^d u_i^{n_i}=0$ locally at $x$ ($n_i$'s are natural numbers). More generally, we say
that $Z$ is a {\em monomial divisor} if etale-locally it is a strictly monomial divisor.
\end{defin}

Note that if each exponent $n_i$ is either $1$ or $0$ in the above definition, then one obtains the usual
definition of (strictly) normal crossing divisor. So, a reduced (strictly) monomial divisor is the same as a
(strictly) normal crossing divisor. Furthermore, the following result holds.

\begin{lem}
\label{snclem} Assume that $X$ is a regular scheme and $Z$ is a closed subscheme, then $Z$ is a (strictly)
monomial divisor if and only if $Z$ is a Cartier divisor and the reduction of $Z$ is a (strictly) normal
crossing divisor in $X$.
\end{lem}
\begin{proof}
Assume that $Z$ is a Cartier divisor at $x$ given by an element $f\in\calO_{X,x}$. Since the regular ring
$\calO_{X,x}$ is factorial, $f=\prod_{i=1}^d f_i^{n_i}$ where each $f_i$ defines an irreducible component of the
reduction of $Z$. Thus the reduction of $Z$ is defined by $\prod_{i=1}^d f_i$ and then it is obvious that $Z$ is
(strictly) monomial if and only if its reduction is (strictly) normal crossing.
\end{proof}

\begin{defin}\label{reglocdef}
(i) Let $X$ and $Z$ be as in the previous definition. The {\em regular locus} $(X,Z)_\reg$ of the pair $(X,Z)$
is the set of points $x\in X$ such that $\calO_{X,x}$ is a regular ring and $Z\times_X\Spec(\calO_{X,x})$ is a
monomial divisor. The {\em singular locus} of the pair $(X,Z)$ is defined as $(X,Z)_\sing=X\setminus(X,Z)_\reg$.

(ii) The {\em strictly regular locus} $(X,Z)_\sreg$ is defined similarly, but with $Z\times_X\Spec(\calO_{X,x})$
being strictly monomial; its complement will be denoted $(X,Z)_\ssing$ (one could call it the semi-singular
locus, but we prefer not to multiply entities beyond necessity).

(iii) If $Z$ is a closed subset, then we set $(X,Z)_\sing=(X,\calZ)_\sing$ and $(X,Z)_\ssing=(X,\calZ)_\ssing$,
where $\calZ$ is the reduced closed subscheme corresponding to $Z$.
\end{defin}

We remark that if $Z'\subseteq Z$ are two Cartier divisors in $X$ then $(X,Z')_\sing\subseteq(X,Z)_\sing$, but
such inclusion is false for general closed subschemes or (even) subsets $Z'\subset Z$.

\begin{lem}\label{openlem}
If $X$ is a quasi-excellent scheme with a closed subscheme $Z$, then $(X,Z)_\reg$ is open. If, furthermore, $X$
is integral and $Z \neq X$, then $(X,Z)_\sreg$ is non-empty.
\end{lem}
\begin{proof}
We assume that $X$ is integral and $Z\neq X$, the general case is proved similarly. We may replace $X$ with any
neighborhood $X'$ of $(X,Z)_\reg$ because $(X',Z\times_X X')_\reg=(X,Z)_\reg$ (we identify $X'$ with a subset of
$X$). So, first of all, replace $X$ with $X_\reg$. Next, removing from $X$ all embedded components of $Z$ and
irreducible components of $Z$ of codimension larger than $1$ (these components lie in $(X,Z)_\sing$), we achieve
that $Z$ is a Cartier divisor. By the previous lemma, we can now replace $Z$ with its reduction.

It now suffices to show that if $Z$ is normal crossing at $x$, then it is normal crossing in a neighborhood of
$x$. Let us first assume that $Z$ is snc at $x$ (as usually, snc stands for strictly normal crossing), and
$Z_1\. Z_m$ are the irreducible components of $Z$ containing $x$. Then, the scheme-theoretic intersection
$T=\cap_{i=1}^mZ_i$ is regular at $x$ and has codimension $m$, hence $x$ possesses a neighborhood $U$ such that
$Z_1\. Z_m$ are the only components of $Z$ which intersect $U$ and $T\cap U$ is regular and has codimension $m$.
It is well known that $Z\cap U$ is snc then. If $Z$ is a normal crossing at $x$, then there exists an etale
neighborhood $f:U\to X$ of $x$ such that $f^{-1}(Z)$ is snc. Since the etale morphism $f$ is an open map, $f(U)$
is a neighborhood of $x$ such that $Z\cap f(U)$ is normal crossing.
\end{proof}

We see that if $X$ is quasi-excellent, then $(X,Z)_\sing$ is a closed subset. Sometimes it will be convenient to
consider it as a reduced closed subscheme.

\begin{lem}
\label{regmorlem} Let $X$ be a quasi-excellent scheme with a closed subscheme $Z$ and $f:X'\to X$ be a regular
morphism. Then $(X',Z')_\sing\toisom(X,Z)_\sing\times_X X'$, where $Z'=Z\times_X X'$.
\end{lem}
\begin{proof}
Since $f$ is regular and the singular loci are reduced, it suffices to check that
$f^{-1}((X,Z)_\sing)=(X',Z')_\sing$ set-theoretically. By \cite[23.7]{Mat}, $f^{-1}(X_\sing)=X'_\sing$, hence we
can replace $X$ with $X_\reg$ and shrink the other schemes accordingly. Obviously, if $T$ is either an embedded
component of $Z$ or an irreducible component of codimension larger than one, then $T'=f^{-1}(T)$ is an analogous
component of $Z'$, hence $T\subset(X,Z)_\sing$ and $T'\subset(X',Z')_\sing$. So, we can remove all such
components from $X$, and then $Z$ becomes a Cartier divisor and using lemma \ref{snclem} we can replace $Z$ with
its reduction.

Let $x\in Z$ be a point and $x'\in Z'$ be its preimage. It suffices to prove that $Z$ is normal crossing at $x$
if and only if $Z'$ is normal crossing at $x'$. Find an etale neighborhood $g:U\to X$ of $x$ such that $u\in U$
is the only preimage of $x$ and all irreducible components of $Z_U=g^{-1}(Z)$ are unibranch at $u$. Then $Z_U$
is normal crossing at $u$ if and only if it is strictly so, hence $Z_U$ is snc at $u$ if and only if $X$ is
normal crossing at $x$. Since the induced morphism $Z'_U=Z_U\times_X X'\to Z_U$ is regular and $Z_U$ is snc at
$u$ if and only if the scheme-theoretic intersection of relevant irreducible components of $Z$ is regular and of
correct codimension, we obtain that $Z_U$ is snc at $u$ if and only if $Z'_U$ is snc at the preimage of $u$
which sits over $x'$. It follows that $Z$ is normal crossing at $x$ if and only if $Z'$ is so at $x'$, as
stated.
\end{proof}

Note that the above lemma fails for strictly regular loci (and their complements $(X,Z)_\ssing$). This fact
forces us to define desingularizations using arbitrary monomial divisors instead of strictly monomial ones. Next
we introduce a notion of desingularization.

\begin{defin}\label{desdef}
(i) Given a locally noetherian scheme $X$ with a closed subscheme $Z$ and a blow up $f:X'\to X$ with support in
$X_\sing\cup Z_\sing$, we say that {\em $f$ desingularizes} the pair $(X,Z)$ {\em over} a subset $S\subset X$ if
$f^{-1}(S)\subset(X',Z\times_X X')_\reg$. If $S=X$ (resp. $S=X^{<d}$), then we say that $f$ {\em desingularizes}
the pair $(X,Z)$ (resp. {\em up to codimension $<d$}), or that $f$ is a {\em desingularization} of the pair. By
a {\em desingularization} of $X$ we mean a desingularization of the pair $(X,\emptyset)$.

(ii) If $Z\subset X$ is a closed subset, then by a desingularization of the pair $(X,Z)$ we mean an
$(X,Z)_\reg$-admissible blow up $f:X'\to X$ such that $(X',f^{-1}(Z))_\reg=X'$. We define desingularizations up
to codimension $<d$ similarly.
\end{defin}

For example, if $(X,Z)_\reg$ is empty (e.g. $Z=X$, or $X$ is not reduced at its maximal points), then
$\Bl_X(X)=\emptyset$ is a desingularization of the pair $(X,Z)$.

\begin{rem}\label{desrem0}
We require $f$ to be a blow up for the following reason. On one hand, blow ups form a sufficiently generic class
of modifications, including, for example, any projective modification of a schemes with an ample sheaf. In
particular, usually it does not cost a serious extra-work to achieve this extra-condition. On the other hand,
working with blow ups one enjoys even more flexibility than when working with a wider class of arbitrary
projective modifications. For example, any blow up of an open subscheme $U\into X$ easily extends to a blow up
of $X$ (while to extend a general modification one would have to invoke Nagata compactification theorem).
\end{rem}

\begin{rem}\label{desrem1}
Usually, by embedded desingularization of $Z$ in $X$ one understands a desingularization of a pair $(X,Z)$ with
a regular ambient scheme $X$. However, it is convenient for our purposes to extend this notion to arbitrary
pairs $(X,Z)$ (for example, we will make some use of reducible $X$'s). In classical terminology,
desingularization of such a pair can be obtained from a desingularization of $X$ and a subsequent embedded
desingularization of the preimage of $Z$. Note also that unlike the case of algebraic varieties, there exist
singular quasi-excellent schemes that even locally cannot be embedded into regular ones.
\end{rem}

\begin{rem}\label{desrem2}
Note that desingularization in our sense provides a control on the exceptional set. For example, a
desingularization $f:X'\to X$ of a pair $(X,X_\sing)$ provides a desingularization of $X$ whose exceptional set
$E$ is a normal crossing divisor (recall that $E$ is the minimal closed set such that $f$ is an open immersion
on $X'\setminus E$, hence $E=f^{-1}(X_\sing)$ in our situation). On the other hand, by a weak desingularization
one usually means a modification $f:X'\to X$ such that $(X',f^{-1}(Z))_\sing=\emptyset$. Perhaps the lack of
control on the exceptional set is the main weakness of weak desingularization.
\end{rem}

Since, usually one imposes more restrictive conditions in the definition of a desingularization, we suggest the
following terminology. See, also remark \ref{resrem}, for the definition of functorial resolution of
singularities.

\begin{defin}
Let a blow up $f:X'\to X$ be a desingularization of a pair $(X,Z)$, then we say that

(i) $f$ is {\em successive} if it is a composition of blow ups along regular centers.

(ii) $f$ is {\em strict} if it is $(X,Z)_\sing$-supported;

(iii) $f$ is {\em semi-strict} if it is $(X,Z)_\ssing$-supported;
\end{defin}

Note that the (semi-) strictness condition is essential only in the embedded case (i.e. when $Z$ is non-empty).
Usually, this extra-condition is not established in desingularization theorems, but it looks a very natural
condition in view of our definition \ref{desdef}. The author is indebted to E. Bierstone and P. Milman for
communicating the proof of the following theorem, which establishes functorial successive strict embedded
desingularization.

\begin{theor}[Bierstone-Milman]\label{BMth}
Let $k$ be a field of characteristic zero. Then for any $k$-smooth scheme $X$ with a reduced closed subscheme
$Z$ there exists a canonical succession $X_n\to X_{n-1}\to\dots\to X_0=X$ of blow ups along smooth centers
$C_i\into X_i$ such that each $C_i$ lies over $(X,Z)_\sing$ and one has $(X_n,Z\times_X X_n)_\sing=\emptyset$.
\end{theor}
\begin{proof}
Let $Z(k)$ denote the set of points $x\in Z$ such that the formal completion of $Z$ along $x$ contains exactly
$k$ irreducible components. Then for a point $x\in Z(k)$ the following properties are easily verified: (i) $x$
possesses a neighborhood $U$ such that the set of points $y\in U$ with $\ord_y Z\ge\ord_x Z$ is contained in
$U\cap Z(k)$, and (ii) $\inv_Z(x)=(k,0,1,0\. 0,1,0,\infty)$ with exactly $k-1$ pairs $1,0$ if and only if $Z$ is
a normal crossing divisor at $x$ (where $\inv$ is the desingularization invariant from \cite{BM}).

Given a sequence of blow ups $X_i\to X_{i-1}\to\dots\to X_0=X$, we denote the composite blow up by $f_i:X_i\to
X$ and define $Z_i\subset X_i$ to be the strict transform of $Z$. If $r$ is the maximal number such that the set
$Z(r)\cap (X,Z)_\sing$ is not empty, then applying the desingularization algorithm of \cite{BM} to the
hypersurface corresponding to $Z$, we can successively blow up along smooth centers with
$\inv_Z(\cdot)>(r,0,1,0\. 1,0,\infty)$ (where the number of pairs $(1,0)$ is $r-1$), until either
$Z_n(r)=\emptyset$ or $\inv_Z(\cdot)=(r,0,1,0\. 1,0,\infty)$ on $Z_n(r)$ for some $n$. In the latter case,
$f_n^{-1}(Z)$ is normal crossing at all points of $Z_n(r)$, hence in a neighborhood of $Z_n(r)$. All in all,
$(X_n,Z_n)_\sing$ is disjoint from $Z_n(r)$, and iterating the same process for smaller values of $r$ we can
achieve that $(X_n,Z_n)_\sing$ is disjoint from $Z_n$. It follows that $f_n^{-1}(Z_n)$ is normal crossing, as
required, and by the construction all centers of successive blow ups lie over $(X,Z)_\sing$. So, the total blow
up $f_n$ is $(X,Z)_\sing$-supported by lemma \ref{compblow}.
\end{proof}

\subsection{Resolution of singularities over a scheme}\label{ressec}

\begin{defin}
\label{resdef} Let $k$ be a locally noetherian scheme. We say that {\em there is resolution of singularities
over $k$ (up to dimension $<d$)} if any integral $k$-scheme $X$ of finite type admits a desingularization (up to
codimension $<d$). If, moreover, for any closed subscheme $Z$, the pair $(X,Z)$ admits a desingularization
$f:X'\to X$ (up to codimension $<d$), then we say that {\em there is embedded resolution of singularities over
$k$ (up to dimension $<d$)}. We say that resolution of singularities over $k$ is {\em strict} (resp. {\em
semi-strict}) if one can choose $f$ to be $(X,Z)_\sing$-supported (resp. $(X,Z)_\ssing$-supported).
\end{defin}

Note that if there exists resolution of singularities over $k$ up to dimension $<d$ then any $k$-scheme of
finite type of dimension strictly less than $d$ admits a desingularization. We will not need the following
definitions, so we put them into a remark.

\begin{rem}\label{resrem}
One can define successive resolution over $k$  in a similar way. Furthermore, we say that there is {\em
functorial} resolution over $k$ (resp. functorial successive resolution over $k$, etc.) if the resolutions
$f_{(X,Z)}:X'\to X$ can be given in a functorial with respect to smooth (or all regular) morphisms way). For
example, minimal resolution of two-dimensional schemes is functorial but not successive, and modern
desingularization theorems provide functorial resolution of algebraic varieties (see \cite[3.106]{Kol} for an
example regarding successiveness).
\end{rem}

The following statement will often be used implicitly.

\begin{lem}
\label{impliclem} Let $X$ be a locally noetherian integral scheme with a closed subscheme $Z$ and $V=(X,Z)_\reg$
(resp. $V=(X,Z)_\sreg$), $f:X'\to X$ be a $V$-admissible blow up, $Z'=Z\times_X X'$ and $d$ be a number. If
$f':X''\to X'$ strictly (resp. semi-strictly) desingularizes $(X',Z')$ up to codimension $<d$, then $f''=f\circ
f':X''\to X$ strictly (resp. semi-strictly) desingularizes $(X,Z)$ up to codimension $<d$.
\end{lem}
\begin{proof}
The two cases are proved in the same way, so we consider only the strict case. Note that $V\otisom
f^{-1}(V)\subset (X',Z')_\reg$, therefore $f'$ is $V$-admissible. By lemma \ref{compblow}, the composition
$f\circ f'$ is $V$-admissible. Obviously $Z'\times_{X'}X''\toisom Z\times_X X''$, hence
$f'^{-1}(X'^{<d})\subset(X'',Z\times_X X'')_\reg$. It remains to note that $f^{-1}(X^{<d})\subset X'^{<d}$ by
the dimension inequality \cite[5.5.8]{EGAIV}. Indeed, for any point $x'\in X'$ with $x=f(x')$, we have that
$\dim(\calO_{X',x'})\le\dim(\calO_{X,x})$ because $\trdeg_K(L)=0$ for $K=\Frac(\calO_{X,x})$ and
$L=\Frac(\calO_{X',x'})$.
\end{proof}

The following proposition allows us to build a global desingularization from local ones. We say that a local
scheme $S$ is of {\em essentially finite type} over a scheme $k$ if it is $k$-isomorphic to $\Spec(\calO_{X,x})$
for a finite type $k$-scheme $X$ with a point $x\in X$.

\begin{prop}
\label{locdesprop} Let $k$ be a noetherian quasi-excellent scheme and $d$ be either a natural number or
infinity, then the following conditions are equivalent.

(i) For any integral $k$-scheme $X$ of finite type with a closed subset $Z$, the pair $(X,Z)$ admits a strict
desingularization up to codimension $<d$.

(ii) There is strict embedded resolution of singularities over $k$ up to dimension $<d$, i.e. for any integral
$k$-scheme $X$ of finite type with a closed subscheme $Z$, the pair $(X,Z)$ admits a strict desingularization up
to codimension $<d$.

(iii) If $S$ is an integral local $k$-scheme of essentially finite type, $\dim(S)<d$, $s\in S$ is the closed
point, $f:S'\to S$ is a blow up and $Z'\into S'$ is a closed subscheme with $(S',Z')_\sing\subset f^{-1}(s)$,
then the pair $(S',Z')$ admits a strict desingularization.

(iv) If $S$ is an integral local $k$-scheme of essentially finite type, $\dim(S)<d$, $s\in S$ is the closed
point, $f:S'\to S$ is a blow up and $Z'\subset S'$ is a closed subset with $(S',Z')_\sing\subset f^{-1}(s)$,
then the pair $(S',Z')$ admits a strict desingularization.

Similar conditions are equivalent when: (1) the resolution is not embedded and $Z=Z'=\emptyset$, (2) the
resolution is embedded, (3) the resolution is embedded and semi-strict.
\end{prop}

\begin{proof}
We consider only the strict embedded case because it is slightly more involved.

(i)$\Rightarrow$(ii) Let $X$ be an integral $k$-scheme of finite type with a closed subscheme $Z\into X$. The
blow up $\Bl_Z(X)\to X$ is an isomorphism over $V=(X,Z)_\reg$ because $Z\times_X V$ is a Cartier divisor in $V$.
By lemma \ref{admisslem}, $\Bl_Z(X)$ is dominated by a $V$-admissible blow up $X'\to X$. Then $Z'=Z\times_X X'$
is a Cartier divisor in $X'$ because already $Z\times_X\Bl_Z(X)$ is a Cartier divisor in $\Bl_Z(X)$. Let
$g:X''\to X'$ be a strict desingularization of $(X',|Z'|)$ up to codimension $<d$. Since $Z''=Z'\times_{X'} X''$
is a Cartier divisor, $(X'',Z'')_\sing=(X'',|Z''|)_\sing$ by lemma \ref{snclem}. It follows that $g$ strictly
desingularizes $(X',Z')$ up to codimension $<d$, and by lemma \ref{impliclem}, the morphism $X''\to X$ strictly
desingularizes $(X,Z)$ up to codimension $<d$.

(ii)$\Rightarrow$(iii) Let $S,S',Z'$ be as in (iii). Find a $k$-scheme $X$ of finite type with a point $x\in X$
such that $\Spec(\calO_{X,x})\toisom S$, in particular, $S$ is a pro-open pro-subscheme of $X$. Replacing $X$
with a neighborhood of $x$ we can make it integral. Furthermore, it follows from \cite[8.6.3]{EGAIV} that
shrinking $X$ further, we can achieve that $f$ is induced from a blow up $X'\to X$, and $Z'\toisom Y'\times_{X'}
S'$ for a closed subscheme $Y'\into X'$. Then $S'\toisom S\times_X X'$ can be identified with a pro-open
pro-subscheme of $X'$ and $Z'=Y'\cap S'$ as sets. For any point $x'\in S'$ we have that
$\dim(\calO_{X',x'})=\dim(\calO_{S',x'})\le\dim(S')$ and $\dim(S')\le\dim(S)$ by the dimension inequality,
\cite[5.5.8]{EGAIV}. Since $\dim(S)<d$, we obtain that $S'\subset (X')^{<d}$. By (ii), the pair $(X',Y')$ can be
strictly desingularized up to codimension $<d$ by a blow up $g:X''\to X'$. Then the pro-open pro-subscheme
$S''=S'\times_{X'}X''$ of $X''$ is regular and its closed subscheme
$Z''=Z'\times_{S'}S''\toisom(Y'\times_{X'}X'')\times_{X''}S''$ is a monomial divisor. Since $g$ is an
$(X',Y')_\reg$-admissible blow up and $(S',Z')_\reg=(X',Y')_\reg\cap S'$ as sets, the morphism $S''\to S'$ is an
$(S',Z')_\reg$-admissible blow up. Therefore, $(S'',Z'')$ is a strict desingularization of $(S',Z')$, and we
obtain (ii).

(iv) follows obviously from (iii), so, it remains to establish the implication (iv)$\Rightarrow$(i). Until the
end of the proof we consider only desingularizations of (scheme, subset) pairs. Set $V=(X,Z)_\reg$, and let
$f:X'\to X$ be a $V$-admissible blow up which desingularizes the pair $(X,Z)$ over an open subscheme $U\into X$.
Assume that $U$ does not contain $X^{<d}$. If we were able to prove that there exists a $V$-admissible blow up
$X''\to X$, which desingularizes $X$ over an open subscheme $W$ with $U\varsubsetneq W$, then the statement of
(i) would follow by noetherian induction (we can start the induction with $f=\Id_X$ and $U=V$).

Let $x\in X^{<d}$ be a maximal point of $X\setminus U$. Set $S=\Spec(\calO_{X,x})$, $Z'=f^{-1}(Z)$,
$S'=S\times_X X'$, $Z_S=Z\cap S$, $Z'_S=f^{-1}(Z_S)$. Then the set $T'=(S',Z'_S)_\sing$ equals to
$(X',Z')_\sing\cap S'$ because $S'$ is a pro-open pro-subscheme of $X'$. Note that $T'\subset f^{-1}(x)$ because
$f^{-1}(S\setminus\{x\})\subset f^{-1}(U)\subset(X',Z')_\reg$. Since $\dim(S')\le\dim(S)<d$, we can apply (iv)
to find a morphism $g:S''\to S'$ which strictly desingularizes $(S',Z'_S)$. The scheme $S''$ is obtained from
$S'$ by blowing up an ideal $\calI\subset\calO_{S'}$ supported on $T'$. Applying lemma \ref{extendlem} to $X'$
and $S'$, we can extend $\calI$ to an ideal $\calJ\subset\calO_{X'}$ supported on the Zariski closure of $T'$.
It follows that $f':X''=\Bl_\calJ(X')\to X'$ is a $U$-admissible blow up, which induces the blow up $g:S''\to
S'$. Therefore, $f''=f\circ f'$ is a $V$-admissible blow up which coincides with $f$ over $U$ and desingularizes
$(X,Z)$ over $x$, i.e. $R=f''((X'',Z'')_\sing)$ is disjoint from the set $U\cup\{x\}$. By properness of $f''$,
the set $R$ is closed, hence $W=X\setminus R$ is as required.
\end{proof}

Combining the implication (i)$\Ra$(ii) from the proposition with theorem \ref{BMth}, we obtain the following
corollary from the results of \cite{BM}.

\begin{cor}\label{BMcor}
There is strict embedded resolution of singularities over any field of characteristic zero.
\end{cor}

In a more general context, Hironaka proved in \cite{Hir} that there is embedded resolution of singularities over
a local quasi-excellent scheme $k$ of characteristic zero: Main Theorem 1 establishes non-embedded
desingularization for schemes of finite type over $k$, and Corollary 3 to the Main Theorem 2 desingularizes
pairs $(X,Z)$ with a regular $X$ of finite type over $k$. Accordingly to proposition \ref{locdesprop},
Hironaka's result implies the following theorem.

\begin{theor}
\label{Groclaim} Any pair of quasi-excellent schemes $(X,Z)$ of characteristic zero with an integral $X$ admits
a desingularization.
\end{theor}

Note that one could also deduce semi-strict embedded desingularization, but then one would have to apply
proposition \ref{standardprop} from the appendix. We finish the section with the following easy lemma.

\begin{lem}
\label{dummylem} Assume that there is (semi-) strict embedded resolution of singularities over $k$. Then for any
$k$-scheme $X$ of finite type over $k$ with a closed subscheme $Z$, the pair $(X,Z)$ admits a (semi-) strict
desingularization.
\end{lem}
\begin{proof}
As usually, we consider only the strict case. Blowing up the ideal generated by all non-zero nilpotent elements
we achieve that $X$ becomes reduced. If $X$ is a disjoint union of integral schemes, then the lemma is trivial.
Assume this is not the case, and let $X_1$ be an irreducible component of $X$ and $X_2$ be the union of all
other components. Then $S=X_1\times_X X_2$ is supported on $X_\sing$, and the blow up of $X$ along $S$ separates
the preimages of $X_i$'s. Applying induction on the number of irreducible components, we obtain that there
exists an $X_\reg$-admissible blow up $X'\to X$ such that $X'$ is a disjoint union of integral schemes. Then, it
suffices to desingularize the pair $(X',Z\times_X X')$, and we are done.
\end{proof}

\begin{rem}
\label{dummyrem} The lemma shows that desingularization of pairs $(X,Z)$ reduces to the case of an integral $X$.
However, it is a stupid desingularization: it brutally separates irreducible components and kills non-reduced
and embedded ones. If one wants to study not integral $X$'s deeper, then the regular locus $(X,Z)_\reg$ should
be replaced with a finer notion. For example, it seems natural to weaken the notion of regularity so that
$X_{"\reg"}=X$ when $X$ is either a reduced strictly normal crossing scheme, or an irreducible scheme normally
flat along $X_\red$.
\end{rem}

\section{Desingularization of special formal schemes}
\label{formch}  Our next and main aim is to prove theorem \ref{mainth} using only resolution of singularities
over fields as proved in \cite{BM} (we will refer only to corollary \ref{BMcor} which is based on theorem
\ref{BMth}). Alternatively, one could use any modern desingularization theorem and the result of the appendix.

\subsection{Quasi-excellent formal schemes and regularity conditions.}
\label{formregsec} If $A$ is a quasi-excellent adic ring, then any formal localization homomorphism $\phi_f:A\to
A_{\{f\}}$ is regular. Indeed, $\phi_f$ is a composition of the localization homomorphism $A\to A_f$ and the
completion homomorphism $A_f\to A_{\{f\}}$, but regularity is preserved by compositions, the first homomorphism
is obviously regular and the second one is regular by \cite[7.8.3(v)]{EGAIV} because $A_f$ is quasi-excellent.
Since there is no published proof in the literature of the fact that the rings $A_{\{f\}}$ are quasi-excellent,
we are forced to give the following definition: a formal scheme $\gtX$ is called {\em absolutely (quasi)
excellent} if for any open affine subscheme $\Spf(A)\into\gtX$ the ring $A$ is (quasi) excellent. We say that
$\gtX$ is {\em (quasi) excellent} if it admits an open covering by absolutely (quasi) excellent subschemes.

\begin{rem}\label{Gabrem}
(i) Since the notion of excellent schemes was defined in \cite[\S7]{EGAIV}, it was an important open question
whether (quasi-) excellence of $A$ implies that of $A[[T]]$ (see loc.cit. 7.8.1.A). In particular, it was not
clear whether completed localization $A_{\{f\}}$ must inherit (quasi-) excellence.

(ii) The author thanks the referee for informing him that recently a (much stronger) ultimate result on
quasi-excellence of adic rings was proved by Offer Gabber: a noetherian complete $I$-adic ring $A$ is
quasi-excellent if and only if $A/I$ is so.

(iii) Gabber's result implies that any quasi-excellent formal scheme is absolutely and universally (see
\S\ref{mainsec}) so. Unfortunately, no printed proof is currently available, so we prefer to make use of these
(superfluous) notions in our paper.
\end{rem}

Many properties of quasi-excellent formal schemes can be defined and studied via their scheme analogs (compare
to \cite[\S1.2]{Con2}; see also \cite[\S2.2]{Ber1}, where one similarly studies $k$-analytic spaces). Consider
an absolutely quasi-excellent affine formal scheme $\gtX=\Spf(A)$ and a closed formal subscheme $\gtZ$ given by
an ideal $\gtI$, and set $X=\Spec(A)$ and $Z=\Spec(A/\gtI)$. We define the {\em singular locus}
$(\gtX,\gtZ)_\sing$ of the pair to be the closed formal subscheme attached to the ideal $\gtJ\subset A$ which
defines the closed subscheme $(X,Z)_\sing\into X$. As we noted in remark \ref{adicrem}, such singular locus is
more informative than the subset of $\gtX_s$ it defines. The following lemma shows that singular loci are
compatible with formal localizations of absolutely quasi-excellent formal schemes.

\begin{lem}
\label{loclem} If $f\in A$ is an element, $\gtX_{\{f\}}=\Spf(A_{\{f\}})$ and
$\gtZ_{\{f\}}=\gtZ\times_\gtX\gtX_{\{f\}}$, then
$(\gtX_{\{f\}},\gtZ_{\{f\}})_\sing\toisom(\gtX,\gtZ)_\sing\times_\gtX \gtX_{\{f\}}$.
\end{lem}
\begin{proof}
The homomorphism $A\to A_{\{f\}}$ is regular because $A$ is quasi-excellent. It remains to note that
$\gtZ_{\{f\}}=\Spf(A/\gtI\otimes_A A_{\{f\}})$, hence $(\gtX_{\{f\}},\gtZ_{\{f\}})_\sing=\Spf(A/\gtJ\otimes_A
A_{\{f\}})$ by lemma \ref{regmorlem}.
\end{proof}

The lemma allows to globalize the definition of singular locus of a formal pair $(\gtX,\gtZ)$ with
quasi-excellent $\gtX$. Indeed, for any open affine absolutely quasi-excellent $\gtX'\into\gtX$ with
$\gtZ'=\gtZ\times_\gtX\gtX'$ we defined a closed subscheme $(\gtX',\gtZ')_\sing\into\gtX'$, and the lemma
implies that these subschemes agree on the intersections of affine open formal subschemes. So, the local
singular loci glue to a single closed formal subscheme $(\gtX,\gtZ)_\sing\into\gtX$. We say that $\gtX$ is {\em
regular} (resp. {\em rig-regular}) if $\gtX_\sing$ is empty (resp. a closed subscheme of $\gtX_s$). We say that
$\gtZ$ is a {\em monomial divisor} (resp. a {\em rig-monomial divisor}) if the formal scheme
$\gtZ\times_\gtX(\gtX,\gtZ)_\sing$ is empty (resp. a closed subscheme of $\gtX_s$). In the following remark we
compare the definition of singular locus to the definition from \cite[\S1.2]{Con2}, and take $\gtZ=0$ for
simplicity. We will not use the following side remark (so we skip the argument), but we hope that it might be
instructive for the reader.

\begin{rem}\label{genrem}
(i) The underlying set of $\gtX_\sing$ coincides with the set of points $x\in\gtX$ with not regular local ring
$\calO_{\gtX,x}$. Thus $\gtX$ is regular in our sense if and only if it is regular as a locally ringed space.

(ii) If $\gtX$ is as in \cite[1.2.1]{Con2} then its singular locus in loc.cit. is defined to be the set
$|\gtX_\sing|$. As we saw, this definition leads to the same notion of regularity, but such singular locus is
less informative. For example, it is always contained in $\gtX_s$ by its definition, so it cannot be used to
define rig-regularity. On the other hand, the support of the closed subscheme $\gtX_\sing$ does not have to be
contained in that of $\gtX_s$, so rig-regularity is a non-trivial condition.

(iii) We noted in remark \ref{adicrem} that $\gtX_\sing$ is determined set-theoretically by
$|(\gtX_\sing)^\ad|$. It seems very probable (and can be proved for $\gtX$ of finite type over a DVR) that the
latter set coincides with the singularity locus of $\gtX^\ad$ viewed as a locally ringed space
$(\gtX^\ad,\calO_{\gtX^\ad})$.

(iv) At least for $\gtX$ of finite type over a DVR, rig-regularity means that a formal scheme $\gtX$ is regular
outside of its closed fiber $\gtX_s$ in the sense that its (rigid, analytic or adic) generic fiber $\gtX_\eta$
is regular.
\end{rem}

The following lemma shows that singular loci are compatible with completions.

\begin{lem}
\label{compllem} Let $X$ be a quasi-excellent scheme with closed subschemes $Z$ and $T=(X,Z)_\sing$,
$\calP\subset\calO_X$ be an ideal and $\gtX$, $\gtZ$, $\gtT$ be the $\calP$-adic completions of $X$, $Z$, $T$.
Assume that $\gtX$ is quasi-excellent, then $\gtT\toisom(\gtX,\gtZ)_\sing$.
\end{lem}
\begin{proof}
Since completions and singular loci are compatible with (formal) localizations, the claim of the lemma is local
on $X$. So, we can assume that $X=\Spec(A)$ is affine with an absolutely quasi-excellent completion $\gtX$,
$\calP$ corresponds to an ideal $P\subset A$, $Z=\Spec(A/I)$ and $T=\Spec(A/J)$ for ideals $I,J\subset A$.
Obviously, $\gtX=\Spf(\hatA)$, where $\hatA$ is the $P$-adic completion of $A$. Since $A$ is quasi-excellent,
the homomorphism $A\to\hatA$ is regular by \cite[7.8.3(v)]{EGAIV}. By lemma \ref{regmorlem},
$(\Spec(\hatA),\Spec(\hatA/I\hatA))_\sing$ is isomorphic to $\Spec(\hatA/J\hatA)$, hence we obtain that
$(\gtX,\gtZ)_\sing\toisom\Spf(\hatA/J\hatA)\toisom\gtT$.
\end{proof}

The lemma provides an easy way to construct examples of singular loci of formal schemes. For example, if we take
a scheme $X$ and complete it along $Y\into X_\reg$, then we obtain a regular formal scheme $\gtX=\hatX_Y$. On
the other hand, if $Y\supset X_\sing$ then $\gtX$ is rig-regular, and if $X_\sing$ contains an irreducible
component $Z$ with $Z\subsetneq Y$ and $Y\cap Z\neq\emptyset$, then $\gtX$ is not even rig-regular.

\begin{cor}
\label{complcor} Keep the notation of the lemma and consider the closed subscheme $Y=\Spec(\calO_X/\calP)$ of
$X$, then:

(i) $\gtX$ is regular and $\gtZ$ is a monomial divisor if and only if there exists a regular neighborhood $U$ of
$Y$ such that $Z\times_X U$ is a monomial divisor;

(ii) $\gtX$ is rig-regular if and only if there exists a neighborhood $U$ of $Y$ such that $U\setminus Y$ is
regular.
\end{cor}
\begin{proof}
Clearly $\gtX$ is regular and $\gtZ$ is a monomial divisor iff $(\gtX,\gtZ)_\sing=\emptyset$, and the latter is
equivalent to $(X,Z)_\sing\cap Y=\emptyset$ by the lemma. Since $(X,Z)_\sing$ is closed, the last equality holds
if and only if $(U,Z\times_X U)_\sing=\emptyset$ for a neighborhood $U$ of $Y$. This proves (i). To prove (ii)
we use the following chain of equivalences: $\gtX$ is rig-regular iff $\gtX_\sing$ is a subscheme in $\gtX_s$
iff each irreducible component of $X_\sing$ is either contained in $Y$ or disjoint from $Y$ iff $U_\sing\subset
Y$ for a neighborhood $U$ of $Y$ iff $U\setminus Y$ is regular.
\end{proof}

Similarly to regularity, one can define reducedness of quasi-excellent formal schemes and prove analogs of
lemmas \ref{loclem} and \ref{compllem}. We will need reducedness later. A further generalization, which will not
be used, is given in the remark.

\begin{rem}
\label{bfPrem} Similarly to regularity, if $\gtX$ is a quasi-excellent formal scheme and $\bfP$ is one of the
following standard properties: $R_n$, CI, Gor, $S_n$ (or their combinations like Reg, Nor, CM, Red), then one
can define the non-$\bfP$ locus $\gtX_{\rm non-\bfP}$ as a closed subscheme of $\gtX$. These loci satisfy
analogs of lemmas \ref{loclem} and \ref{compllem}, and everything noted in remark \ref{genrem} is valid for
them.
\end{rem}

\subsection{Special formal schemes}
\label{specsec} Throughout this section $k$ denotes a field with $p=\cha(k)$, and $\gtX$ is a formal scheme such
that $\gtX_s$ is a $k$-scheme. We say that $\gtX$ is {\em equicharacteristic} if $p\calO_\gtX=0$. If $p$ is
positive, then a discrete valuation ring $\Kcirc$ is called a {\em $p$-ring} if $p$ generates the maximal ideal
of $\Kcirc$. It is shown in \cite[Ch. 29]{Mat} that up to an isomorphism, there exists a unique complete
$p$-ring with residue field $k$; it will be denoted $\bfZ_p(k)$.

Let $A$ be a $P$-adically complete noetherian ring and $\gtS=\Spf(A)$. Recall that a formal $\gtS$-scheme $\gtX$
is of {\em finite type} if it admits a finite covering by formal schemes of the form $\Spf(A\{T_1\. T_n\}/I)$.
More generally, by an {\em $A$-special} ring we mean a topological ring $A\{T_1\. T_n\}[[R_1\. R_m]]/I$ provided
with the $Q$-adic topology, where $Q$ is generated by the images of $P$ and $R_1\. R_m$. A formal $\gtS$-scheme
$\gtX$ is called {\em $\gtS$-special} if it admits a finite open covering by formal spectra of $A$-special rings
(a definition in \cite{Ber3} is slightly different). A noetherian formal scheme is called {\em special} if its
closed fiber is a scheme of finite type over a field.

\begin{prop}
\label{specialprop0} Assume that $\gtX$ is an affine special formal scheme such that $\gtX_s$ is of finite type
over $k$.

(i) If $\gtX$ is equicharacteristic, then it is isomorphic to a $k$-special formal scheme.

(ii) If $p>0$, then $\gtX$ is isomorphic to a $\bfZ_p(k)$-special formal scheme.
\end{prop}
\begin{proof}
Let $\gtX=\Spf(B)$ and $P$ be the biggest ideal of definition, so $B/P$ is a finitely generated $k$-algebra. If
$B$ is equicharacteristic (in particular, it is automatically the case when $p=0$), then we set $\Kcirc=k$ and
$\Lcirc=\bfF_p$ (if $p=0$ then $\Lcirc=\bfQ$). Otherwise we set $\Kcirc=\bfZ_p(k)$ and $\Lcirc=\bfZ_p$. We have
natural homomorphisms $\Kcirc\to k\to B/P$ and $\Lcirc\to B$, and the same argument as in the proof of Cohen's
structure theorem given in \cite[29.2]{Mat} shows that there is a lifting $\Kcirc\to B$. Indeed, it is shown in
the proof of loc.cit. that $\Kcirc$ is formally smooth over $\Lcirc$, hence the $\Lcirc$-homomorphism $\Kcirc\to
B/P$ lifts to $B/P^2$, $B/P^3$, etc. Since $B$ is $P$-adically complete, we obtain a lifting $\Kcirc\to B$.

Let $r_1\. r_m$ be generators of $P$ and $t_1\. t_n\in B$ be such that their images in $B/P$ generate it over
$k$. Then there exists a continuous homomorphism $\phi:C=\Kcirc\{T_1\. T_n\}[[R_1\. R_m]]\to B$, which takes
$T_i$ and $R_j$ to $t_i$ and $r_j$. The maximal ideal of definition $Q\subset C$ is generated by $p$ and
$R_j$'s, hence $C/Q\toisom k[T_1\. T_n]$ and the induced homomorphism $C/Q\to B/P$ is surjective. The image of
$Q$ in $B$ generates $P$, hence $B$ is topologically finitely generated over $C$ by \cite[0.7.5.5(a)]{EGAI}.
Moreover, the same argument as in the proof of loc.cit. implies that the homomorphism $C\to B$ is actually
surjective because $C/Q\to B/P$ is onto. So, $\gtX$ is $\Kcirc$-special, as required.
\end{proof}

Any $\Kcirc$-special formal scheme is absolutely excellent by results of Valabrega, see \cite[Prop. 7]{Val1} and
\cite[Th. 9]{Val2}. We obtain the following corollary, which allows to apply results of the previous section to
special formal schemes.

\begin{cor}
Any special formal scheme is excellent.
\end{cor}

The following statement is proved exactly as its analog \ref{specialprop0}.

\begin{prop}
\label{specialprop} Assume that $\gtX=\Spf(B)$, $\gtX_s$ is of finite type over $k$, $B$ possesses a principal
ideal of definition $\pi B$, and either $B$ is equicharacteristic or $\pi=p$. Set $\Kcirc=k$ if $\pi$ is
nilpotent, set $\Kcirc=k[[\pi]]$ if $\pi$ is not nilpotent and $B$ is equicharacteristic, and set
$\Kcirc=\bfZ_p(k)$ otherwise. Then $\gtX$ is isomorphic to a formal $\Kcirc$-scheme of finite type.
\end{prop}

\begin{rem}
\label{mixedrem} There exist special formal schemes with principal ideal of definition not covered by the
proposition. For example, if $B$ equals to $\bfZ_p\{T\}[[R]]/(RT-p)$, then $RB$ is an ideal of definition, but
$pB$ is not.
\end{rem}

\subsection{Rig-smoothness and algebraization in characteristic zero}
\label{rigsmoothsec} Let $\calO$ be a ring. We say that a formal scheme $\gtX$ is {\em $\calO$-algebraizable} if
it is isomorphic to the formal completion of an $\calO$-scheme of finite type along a closed subscheme. We say
that $\gtX$ is {\em locally $\calO$-algebraizable} if it can be covered by open $\calO$-algebraizable formal
subschemes.

Fix the following notation: $K$ is a complete discretely valued field with ring of integers $\Kcirc$, residue
field $k$ and maximal ideal $(\pi)$, $\gtS=\Spf(\Kcirc)$, $\gtX$ is a $\Kcirc$-special formal scheme, $L\subset
K$ is a dense subfield and $\calO=L\cap\Kcirc$.

\begin{prop}
\label{algebprop} Assume that $\cha(K)=0$. If $\gtX$ is rig-regular, affine and of finite type over $\gtS$, then
it is $\calO$-algebraizable.
\end{prop}
\begin{proof}
We have that $\gtX=\Spf(C)$ with $C$ topologically finitely generated over $\Kcirc$. Set $X=\Spec(C)$, then
$X\setminus\gtX_s$ is regular by the definition of rig-regularity, and the $K$-affinoid algebra
$C_K=C\otimes_{\Kcirc} K$ is regular because $\Spec(C_K)=X\setminus\gtX_s$. The regular $K$-affinoid space
$\gtX^\rig=\Sp(C_K)$ is $\Sp(K)$-smooth because $K$ is perfect, see \cite[2.8(b)]{BLR}, or an explanation in the
proof of \cite[3.3.1]{Con2}. It is well known (and will be proved in proposition \ref{jacobprop} below because
of lack of an appropriate reference) that $\gtX^\rig$ is $K$-smooth iff $\gtX$ is formally $\Kcirc$-smooth
outside of $V(\pi)$ in the sense of \cite[Th. 7]{Elk}. Thus, we can apply this theorem of Elkik to $\gtX$ and we
thereby obtain that $C$ is isomorphic to the $\pi$-adic completion of a finitely generated $\calO^h$-algebra
$A$, where $\calO^h$ is the henselization of $\calO$. Since $\calO^h$ is a union of etale $\calO$-algebras, $A$
is isomorphic to an algebra $A'\otimes_{\calO'}\calO^h$, where $\calO'$ is $\calO$-etale and $A'$ is a finitely
generated $\calO'$-algebra. Since the completions of $A$ and $A'$ are canonically isomorphic, we obtain that
$\Spec(A')$ is an $\calO$-algebraization of $\gtX$.
\end{proof}

Let $A$ be a noetherian $P$-adic ring with $P$-adic topologically finitely generated $A$-rings $B=A\{T_1\.
T_n\}$ and $C=B/(f_1\. f_l)$. Let us define a topological Jacobian ideal $H_{C/A}$ similarly to its algebraic
analog from \cite[0.2]{Elk}. Set $\Delta=\left(\frac{\partial f_i}{\partial T_j}\right)_{1\le i\le l,1\le j\le
n}$. For any subset $L\subset\{1\. l\}$ with $|L|=r$, let $H_L\subset B$ be the ideal generated by the
determinants of $r\times r$-minors of $\Delta$ whose rows are numbered by the elements of $L$. Also, let $J_L$
be the ideal generated by $f_i$'s with $i\in L$ and let $J=(f_1\. f_l)$. Then we set
$$H_{C/A}=\sqrt{\sum_{L\subset\{1\. l\}} (J_L:J)H_L}C,$$
where $(J_L:J)=\{x\in B|xJ\subset J_L\}$.

Note that a priori $H_{C/A}$ depends on the choices of $B$ and $f_i$'s. A standard argument using the Jacobian
criterion of smoothness shows that in the algebraic case (i.e. $P$ is nilpotent) $H_{C/A}$ defines the not
$A$-smooth locus of $\Spec(C)$; in particular, it is independent of all choices. We refer to \cite[2.13]{Spi}
for details. A similar argument involving modules $\hatOmega^1_{B/A}\toisom\oplus BdT_j$ and $\hatOmega^1_{C/A}$
of continuous differentials shows that $\Spf(C)$ is $\Spf(A)$-smooth iff $H_{C/A}=C$ (a definition of smooth
morphisms of formal schemes can be found in \cite[1.1]{BL2}). In \cite[Th. 7]{Elk}, $\Spf(C)$ is said to be {\em
formally $A$-smooth outside of $V(\pi)$} if the Jacobian ideal $H_{C/A}$ is open. Using the Jacobian criterion
of rig-smoothness, see \cite[3.5]{BLR}, one can show that it happens iff $\Spf(C)$ is rig-smooth over $A$. For
the sake of simplicity, we consider only the classical rigid case which was used in the previous proposition.
Our proof is an affinoid adjustment of the proof of \cite[2.13]{Spi}.

\begin{prop}
\label{jacobprop} Keep the above notation, assume that $A$ is topologically finitely generated over $\Kcirc$,
and set $\calA=A\otimes_{\Kcirc}K$ and $\calC=C\otimes_{\Kcirc}K$. Then the morphism $\Sp(\calC)\to\Sp(\calA)$
of rigid affinoid spaces is smooth if and only if the Jacobian ideal $H_{C/A}$ is open.
\end{prop}
\begin{proof}
Obviously, $\calB/(f_1\. f_l)\toisom\calC$, where $\calB=\calA\{T_1\.T_n\}$. Use this representation of $\calC$
to define the Jacobian ideal $H_{\calC/\calA}$ of affinoid algebras analogously to its adic analog $H_{C/A}$.
Note that the definition of the Jacobian ideal is compatible with localization by $\pi$, hence
$H_{\calC/\calA}=H_{C/A}\calC$, and we obtain that $H_{C/A}$ is open iff $H_{\calC/\calA}=\calC$.

It remains to prove that $H_{\calC/\calA}$ defines the not $\calA$-smooth locus of $X=\Sp(\calC)$. Recall that
modules of differentials of rigid spaces are defined by use of modules of continuous differentials of affinoid
algebras, see \cite[\S1]{BLR}. For example, $\hatOmega_{\calB/\calA}=\oplus_{i=1}^n\calB dT_i$, though
$\Omega_{\calB/\calA}$ can be huge. Let $J\subset\calB$ be the ideal generated by $f_i$'s, then by
\cite[1.2]{BLR}, there is a natural sequence of finite $\calC$-modules (perhaps not exact on the left)
$$
0\to J/J^2{\stackrel{d_{\calA/\calB/\calC}}{\to}}\hatOmega_{\calB/\calA}\otimes_\calB\calC\to
\hatOmega_{\calC/\calA}\to 0
$$
Set $d=d_{\calA/\calB/\calC}$ for shortness. Let $x\in X$ be a point and $m\subset\calC$ be the corresponding
ideal. By the Jacobian criterion of smoothness, see \cite[2.5]{BLR}, $X$ is $\calA$-smooth at $x$ iff the above
sequence becomes split exact after tensoring with $\calO_{X,x}$, or, that is equivalent, the map
$d\otimes_\calC\calO_{X,x}$ has left inverse. Since $\calO_{X,x}$ is local and $\hatOmega_{\calB/\calA}$ is
free, the latter can happen iff the tensored sequence is an exact sequence of free $\calO_{X,x}$-modules.

Suppose that $X$ is $\calA$-smooth at $x$, then there exists a subset $L\subset\{1\. l\}$ such that: (i)
$J/J^2\otimes_\calB\calO_{X,x}$ is freely generated by the images of the elements of $f_L=\{f_i\}_{i\in L}$, and
(ii) the elements of $df_L$ are linearly independent modulo $m$. Identify $X$ with the closed subspace of
$Y=\Sp(\calB)$, then (i) implies that the image of $f_L$ generates $J/J^2\otimes_\calB K(x)\toisom
J\calO_{Y,x}/m_{Y,x}J\calO_{Y,x}$, where $K(x)=\calO_{Y,x}/m_{Y,x}$ is the residue field of $x$. Therefore,
$f_L$ generates the $\calO_{Y,x}$-module $J\calO_{Y,x}$ by lemma of Nakayama, i.e.
$J_L\calO_{Y,x}=J\calO_{Y,x}$. Note that the operation $:$ is compatible with flat base changes (use that
$(I:J)=\Ann(J/I)$), in particular, $(J_L:J)\calD=(J_L\calD:J\calD)$ for any flat $\calB$-algebra $\calD$. Thus,
$(J_L:J)\calO_{Y,x}=(J_L\calO_{Y,x}:J\calO_{Y,x})=\calO_{Y,x}$, and we obtain that $x$ is not contained in
$V((J_L:J))\subset\Spec(\calB)$ (recall that set-theoretically $Y$ coincides with the set of closed points of
$\Spec(\calB)$). Since $df_i=\sum_j\frac{\partial f_i}{\partial T_j}dT_j$, (ii) implies that the rank of
$\left(\frac{\partial f_i}{\partial T_j}(x)\right)_{i\in L,1\le j\le n}$ equals to $|L|$. It follows that
$x\notin V(H_L)$, hence $x\notin V((J_L:J)H_L)$, and, finally, $x\notin V(H_{\calC/\calA})$.

Conversely, suppose that $x\notin V(H_{\calC/\calA})$. Then there exists $L\subset\{1\. l\}$ such that $x\notin
V((J_L:J)H_L\calC)$. Therefore $(J_L:J)H_L\calO_{X,x}=\calO_{X,x}$, and we obtain that
$(J_L:J)H_L\calO_{Y,x}=\calO_{Y,x}$ because $\calO_{X,x}$ is a quotient of the local ring $\calO_{Y,x}$. Then
the set $f_L$ generates $J\calO_{Y,x}$ because
$\calO_{Y,x}=(J_L:J)\calO_{Y,x}\subseteq(J_L\calO_{Y,x}:J\calO_{Y,x})$. Hence the images of the elements of
$f_L$ generate $J/J^2\otimes_\calC\calO_{X,x}$, and, moreover, they generate it freely because their images in
$\hatOmega_{\calB/\calA}\otimes_\calB\calC/m\toisom\oplus_{i=1}^n K(x)dT_i$ are linearly independent (by the
assumption on $H_L$). It follows that $d\otimes_\calC\calO_{X,x}$ has left inverse, hence $X$ is $\calA$-smooth
at $x$.
\end{proof}

Although we do not use that in the sequel, we remark that the Jacobian ideals depend only on the corresponding
homomorphism $A\to C$ or $\calA\to\calC$. Indeed, a reduced closed subspace $Z\subset X$ is uniquely defined by
the set of its points, hence it follows from the above proof that the ideal $H_{\calC/\calA}$ depends only on
the $\calA$-affinoid algebra $\calC$. Moreover, since a reduced closed formal subscheme $\gtZ\subset\Spf(C)$ is
uniquely defined by the sets $\gtZ_s$ and $\gtZ_\eta^\rig$, the Jacobian ideal $H_{C/A}$ depends only on the
homomorphism $A\to C$.

\begin{question}
\label{algebrquest} Assume that $\gtX$ is special rig-regular and admits a locally principal ideal of
definition, and set $\gtT=\gtX_\sing$. Set $\calO=k[\pi]$ if $\gtX$ is equicharacteristic, and let $\calO$ be a
$p$-ring with residue field $k$ otherwise. Does there exist a $\gtT$-supported blow up $\gtX'\to\gtX$ with
locally $\calO$-algebraizable $\gtX'$?
\end{question}

The positive answer to the above question would allow to reduce desingularization of an arbitrary
quasi-excellent scheme $X$ of characteristic $p$ (resp. mixed characteristic) to the particular case of
$k(x)$-schemes of finite type for points $x\in X$ (resp. $\calO$-schemes of finite type, where $\calO$ is a
$p$-ring with residue field $k(x)$).

\subsection{Formal desingularization and applications to schemes}
\label{formmainsec} Given a quasi-excellent formal scheme $\gtX$ with a closed formal subscheme $\gtZ$ and
$\gtT=(\gtX,\gtZ)_\sing$, we say that a $\gtT$-supported blow up $f:\gtX'\to\gtX$ {\em strictly desingularizes}
the pair $(\gtX,\gtZ)$ {\em over} a subset $S\subset\gtX$ if $f^{-1}(S)$ is disjoint from the underlying
topological space of $(\gtX',\gtZ')_\sing$, where $\gtZ'=\gtZ\times_\gtX\gtX'$. If $S=\gtX$, then we say that
$f$ is a {\em strict desingularization} of the pair, it happens iff $\gtX'$ is regular and $\gtZ'$ is a monomial
divisor. In the following theorem we prove that certain special formal schemes of characteristic zero admit a
desingularization. The theorem will be used in the proof of a more general theorem \ref{formmainth}.

\begin{theor}
\label{formtheor} Let $\gtX$ be a reduced rig-regular special formal scheme of characteristic zero with a
locally principal non-zero ideal of definition and $\gtZ$ be a closed $\gtX_s$-supported subscheme. Then the
pair $(\gtX,\gtZ)$ admits a strict desingularization.
\end{theor}
\begin{proof}
Note that $\gtT$ is a reduced closed subscheme of $\gtX_s$, so we can and will identify it with a closed subset
of $\gtX_s$. For any $\gtT$-supported formal blow up $\hatf:\gtX'\to\gtX$, the singular locus
$\gtT'=(\gtX',\gtZ')_\sing$ is $\gtT$-supported (it suffices to check this claim locally on $\gtX$, but if
$\gtX$ is affine, then the statement follows from its analog for schemes). Thus, we can identify $\gtT'$ with a
closed subset of $\gtX'_s$. Assume that $\hatf$ strictly desingularizes the pair $(\gtX,\gtZ)$ over an open
formal subscheme $\gtU$ with $\gtX\setminus\gtT\subseteq\gtU\varsubsetneq\gtX$, for example $\hatf=\Id_\gtX$ and
$\gtU=\gtX\setminus\gtT$. By noetherian induction, it suffices to prove that there exists a $\gtT$-supported
formal blow up $\gtX''\to\gtX$, which desingularizes $\gtX$ over an open formal subscheme $\gtW$ with
$\gtU\varsubsetneq\gtW$.

Choose a field $k$ such that $\gtX_s$ is of finite type over $k$ and set $\Kcirc=k[[\pi]]$ and $\calO=k[\pi]$.
Find an open affine subscheme $\gtX_0$ which possesses a principal ideal of definition and has a non-empty
intersection with the set $\gtS=\gtX\setminus\gtU$. Set $\gtX'_0=\gtX_0\times_\gtX\gtX'$,
$\gtZ_0=\gtZ\times_\gtX\gtX_0$, $\gtU_0=\gtU\times_\gtX\gtX_0$ and $\gtS_0=\gtS\cap\gtX_0$, and let
$\hatf_0:\gtX'_0\to\gtX_0$ be the induced formal blow up. By proposition \ref{specialprop}, $\gtX_0$ is
isomorphic to a formal $\Kcirc$-scheme of finite type, hence $\gtX_0$ is $\calO$-algebraizable by proposition
\ref{algebprop}. Say, $\gtX_0\toisom\hatX_\Pi$ where $X$ is an $\calO$-scheme of finite type and $\Pi\subset X$
is the divisor defined by $\pi$. The closed subscheme $\gtZ_0\into\gtX_0$ is supported on the closed fiber of
$\gtX_0$, so it is given by an open ideal in $\calO_{\gtX_0}$ (in particular, its ideal of definition is
nilpotent, so $\gtZ_0$ is a usual scheme). Therefore $\gtZ_0$ algebraizes to a closed subscheme $Z\into X$
supported on $\Pi$, i.e. $\gtZ_0\toisom\hatZ_\Pi\toisom Z$. By corollary \ref{complcor}, replacing $X$ with a
neighborhood of $\Pi$ we achieve that $X$ is reduced and $X\setminus\Pi$ is regular. Then $T=(X,Z)_\sing$ lies
in $\Pi$ and $T$ is isomorphic to $\gtT_0=(\gtX_0,\gtZ_0)_\sing$ by lemma \ref{compllem}. Let $S\subset\Pi$ be
the preimage of $\gtS_0$ under the homeomorphism $\Pi\to\gtX_0$ and $U=X\setminus S$, then $\gtU_0$ is
isomorphic to the formal completion of $U$ along $U\cap\Pi$.

The formal scheme $\gtX'_0$ is obtained from $\gtX_0$ by blowing up an open ideal $\gtI$ supported on
$\gtT_0=\gtT\cap\gtX_0$, hence $\gtI$ is the completion of a $T$-supported ideal $\calI\subset\calO_X$. If
$f:X'\to X$ denotes the blowing up along $\calI$ and $\Pi'=\Pi\times_X X'$, then $\hatX'_{\Pi'}\toisom\gtX'_0$
and $\hatf_0$ is the completion of $f$ by lemma \ref{formblowlem}. Since $\hatf_0$ strictly desingularizes the
pair $(\gtX_0,\gtZ_0)$ over $\gtU_0$, $f$ strictly desingularizes the pair $(X,Z)$ over $U$ by \ref{complcor}.

Note that $X'$ is reduced because $X$ is. If $X'$ is integral, then the pair $(X',Z\times_X X')$ admits a strict
desingularization $f':X''\to X'$ by corollary \ref{BMcor}. Moreover, by lemma \ref{dummylem}, the integrality
assumption is redundant and $f'$ exists unconditionally. Note that $f'$ is an $S$-supported blow up, hence the
induced morphism $X''\to X$ is a strict desingularization of $(X,Z)$ which is isomorphic to $f$ over $U$.
Passing to completions, we obtain an $\gtS_0$-supported formal blow up $\hatf'_0:\gtX''_0\to\gtX'_0$. The
composition $\gtX''_0\to\gtX_0$ strictly desingularizes $(\gtX_0,\gtZ_0)$ by corollary \ref{complcor} and
coincides with $\hatf_0$ over $\gtU_0$. Using \ref{formextendlem}, we can extend $\hatf'_0$ to an
$\gtS$-supported formal blow up $\hatf':\gtX''\to\gtX'$, in particular, $\hatf'$ is $\gtT$-supported. Then
$\hatf'$ is an isomorphism over $\gtU\times_\gtX\gtX'$ and induces the formal blow up $\gtX''_0\to\gtX'_0$.
Therefore, $\hatf\circ\hatf'$ strictly desingularizes $\gtX$ over $\gtW=\gtU\cup\gtX_0$.
\end{proof}

\begin{cor}
\label{desingcor} Let $X$ be an integral noetherian quasi-excellent scheme of characteristic zero with a closed
subscheme $Z$ such that $X_\sing\subset Z$ and $Z$ is isomorphic to a $k$-scheme of finite type. Then the pair
$(X,Z)$ admits a strict desingularization.
\end{cor}
\begin{proof}
The blow up $\Bl_Z(X)\to X$ is an isomorphism over $V=(X,Z)_\reg$ because $Z\times_X V$ is a Cartier divisor in
$V$. By lemma \ref{admisslem}, $\Bl_Z(X)$ is dominated by a $V$-admissible blow up $X'\to X$, and then
$Z'=Z\times_X X'$ is a Cartier divisor in $X'$. Replacing $X$ and $Z$ with $X'$ and $Z'$, we achieve that $Z$ is
a Cartier divisor.

Let $\gtX$ be the formal completion of $X$ along $Z$, it is reduced and rig-regular by corollary \ref{complcor}
(i). Thus, $\gtX$ satisfies the assumptions of the theorem. The closed subschemes $Z$ and $T=(X,Z)_\sing$ can be
identified with closed subschemes of $\gtX$ because they are supported on $Z_\red\toisom\gtX_s$. Then
$T=(\gtX,Z)_\sing$ by lemma \ref{compllem}.

By the previous theorem, there exists an open ideal $\gtI\subset\calO_\gtX$ supported on $T$ such that
$\gtX'=\Bl_\gtI(\gtX)$ is regular and $Z\times_\gtX\gtX'$ is a monomial divisor. Since $\gtI$ is open, it is the
completion of an ideal $\calI\subset\calO_X$ supported on $T$. Let $X'$ be the blow up of $X$ along $\calI$ and
$Z'=Z\times_X X'$. By lemma \ref{formblowlem}, $\gtX'$ is isomorphic to the formal completion of $X'$ along
$Z'$, hence $Z'$ is a monomial divisor by \ref{complcor} (ii). Since $X'\setminus Z'\toisom X\setminus Z$ is
regular, $X'\to X$ is a required desingularization.
\end{proof}

Now, we obtain a new proof of the not embedded case of theorem \ref{mainth} (it was earlier deduced from
\cite{Hir} in theorem \ref{Groclaim}).

\begin{theor}
\label{desingtheor} Let $k$ be a noetherian scheme of characteristic zero. Then $k$ is quasi-excellent if and
only if there is resolution of singularities over $k$.
\end{theor}
\begin{proof}
The converse implication is due to Grothendieck, see \cite[7.9.5]{EGAIV}. Conversely, by proposition
\ref{locdesprop}, it suffices to prove that if $S$ is an integral local $k$-scheme of essentially finite type,
$s\in S$ is a closed point, and $f:S'\to S$ is a blow up with $S'_\sing\subset f^{-1}(s)$, then $S'$ admits a
desingularization. Note that $S'_\sing$ is of finite type over $k(s)$, hence the pair $(S',S'_\sing)$ admits a
desingularization $g:S''\to S'$ by the previous corollary. Then it is clear that $g$ is a required
desingularization of $S'$.
\end{proof}

\section{Strict transforms and main results}
\label{strch}

The first two sections of \S\ref{strch} are devoted to the proof of proposition \ref{regprop} which
desingularizes strict transforms and is of independent interest. In particular, the proposition is valid for any
ambient scheme $X$ (no restriction on the characteristic) and is used in the appendix. Its proof makes no use of
the results of \S\ref{formch}. Then we combine the proposition with the results of \S\ref{formch} to prove
theorem \ref{mainth}.

\subsection{Principalization of strict transform}
\label{princstrsec}

Let $f:X'=\Bl_T(X)\to X$ be a blow up and $Y$ be a closed subscheme of $X$. We refer to \cite[\S1]{Con1} for an
explicit definition of the {\em strict transform} $\tilY$ of $Y$ in $X'$, but we recall the following property
of $\tilY$ which can be taken as an alternative definition: $\tilY$ coincides with the schematic closure of
$(Y\setminus T)\times_X X'\toisom Y\setminus T$ in $X'$ by \cite[1.1]{Con1} (in particular, the schematic
closure exists). Furthermore, $\tilY$ is canonically isomorphic to the blow up of $Y$ along $T\times_X Y$, see
\cite[\S1]{Con1} for details. If $Y=\Spec(\calO_X/\calI)$ and $Z=\Spec(\calO_X/\calJ)$ are two closed subschemes
and $T=Y\times_X Z\toisom\Spec(\calO_X/(\calI+\calJ))$ is their intersection, then it is well known that blowing
up $X$ along $T$ separates the strict transforms of $Y$ and $Z$, see \cite{Con1}, lemma 1.4 and the consequent
remark. We will also need the following slightly more specific result.

\begin{lem}
\label{strcarlem} Keep the above notation and let $Y'$ be the strict transform of $Y$ in $X'=\Bl_T(X)$. Then
$Z\times_X X'$ is a Cartier divisor in a neighborhood of $Y'$.
\end{lem}
\begin{proof}
Following the proof of \cite[1.4]{Con1}, we assume that $X=\Spec(A)$ and $\calI,\calJ$ correspond to ideals
$I,J\subset A$. It is shown in loc.cit. that $Y'$ is covered by the charts $\Spec(A[\frac{I+J}{f}])$ with $f\in
J$, so it remains to note that the ideal $JA[\frac{I+J}{f}]$ is principal because it coincides with
$fA[\frac{I+J}{f}]$.
\end{proof}

\begin{prop}
\label{princprop} Let $X$ be a noetherian scheme with closed subschemes $T\into Y$ corresponding to
$\calO_X$-ideals $\calJ\supset\calI$. Assume that $Y$ is a Cartier divisor. Given a positive integer $n$, let
$X_n$ denote the blow up of $X$ along $\calJ_n=\calI+\calJ^n$ and let $Y_n$ denote the strict transform of $Y$
in $X_n$. Then the following statements hold true.

(i) $Y_n$ is canonically isomorphic to the blow up of $Y$ along the ideal $\calJ/\calI\subset\calO_Y$.

(ii) $Y_n$ is a Cartier divisor for any sufficiently large $n$.

(iii) Assume that $Y\setminus T$ is a disjoint union of its closed subschemes $\tilY'$ and $\tilY''$, and
$T\toisom Y'\times_X Y''$, where $Y'$ and $Y''$ are the schematic closures of $\tilY'$ and $\tilY''$ in $X$.
Then the strict transforms $Y'_n$ and $Y''_n$ of $Y'$ and $Y''$ in $X_n$ are Cartier divisors for any
sufficiently large $n$.
\end{prop}

\begin{proof}
Recall that $Y_n$ is canonically isomorphic to the blow up of $Y$ along $\calJ_n/\calI$, which equals to the
$n$-th power of $\calJ/\calI$. Since blow ups along an ideal and its powers are canonically isomorphic, we
obtain (i). Let us prove the statement of (ii). It suffices to find a finite covering $\calU$ of $X$ by open
subschemes such that for any $X'\in\calU$ the triple $(X',\calI|_{X'},\calJ|_{X'})$ satisfies (ii). So, we can
assume that $X=\Spec(A)$ is affine and $\calI$ corresponds to a principal ideal $I=fA$. Let $J,J_n\subset A$
denote the ideals corresponding to $\calJ,\calJ_n$. Applying Artin-Rees lemma to the ideal $J$ and the inclusion
of $A$-modules $I\subset A$, we find a positive $n_0$ such that for any $n\ge n_0$ and $m\ge 0$, $J^{n+m}\cap
I=J^m(J^n\cap I)$. Obviously $J^{m+n}\cap I\subset fJ^m$ then. Fix $n\ge n_0$, we will prove that it is as
required. Let us first check that $J_n^k\cap I=fJ_n^{k-1}$ for any $k>0$. Obviously, only the direct inclusion
needs a proof. From the equality $J_n=I+J^n$ we obtain that $J_n^k=fJ_n^{k-1}+J^{nk}$, hence it remains to use
that $J^{nk}\cap I\subset fJ^{n(k-1)}$ (take $m=n(k-1)$ in the above inclusion).

Consider an affine chart $U_g=\Spec(B)$ of $X_n$, where $B=A[\frac{J_n}{g}]$ for some $g\in J_n$. It suffices to
prove that $Y_n\cap U_g$ coincides with the closed subscheme $V_g=V(\frac{f}{g})$ of $U_g$. Note that the
intersection of any of these two schemes with $\Spec(A_g)$ coincides with $Y_g=\Spec(A_g/fA_g)$, and $Y_n\cap
U_g$ is the scheme-theoretical closure of $Y_g$. Therefore, we have only to prove that $Y_g$ is schematically
dense in $V_g$, or, that is equivalent, that the homomorphism $\phi:B/\frac{f}{g}B\to A_g/fA_g$ has no kernel.
Suppose, conversely, that $\phi$ is not injective, then there exists an element $x\in B\setminus\frac{f}{g}B$
such that $x$ is divided by $f$ in $A_g$. Obviously, $x=\sum_{i=0}^l\frac{x_i}{g^i}$ for some elements $x_i\in
J_n^i$, and $g^{l+m}x\in fA$ for sufficiently large $m$. The element $g^{l+m}x=g^m\sum_{i=0}^l x_ig^{l-i}$ is
contained in $J_n^{l+m}$ and is divided by $f$, therefore it is contained in $fJ_n^{l+m-1}$ by the previous
paragraph, i.e. $g^{l+m}x=fh$ for some $h\in J_n^{l+m-1}$. It follows that
$x=\frac{h}{g^{l+m-1}}\frac{f}{g}\in\frac{f}{g}B$, contradicting our assumptions.

It remains to prove (iii). We know from (ii) that $Y_n$ is a Cartier divisor. Note that $Y_n$ is the schematic
closure of $\tilY=Y\setminus T$ in $X_n$. Since $\tilY=\tilY'\sqcup\tilY''$, we obtain that $Y_n=Y'_n\cup
Y''_n$. Therefore, it suffices to prove that $Y'_n$ and $Y''_n$ are disjoint. By part (i), $Y_n$ is isomorphic
to the blow up of $Y$ along $T\toisom Y'\times_X Y''\toisom Y'\times_Y Y''$, but blowing up of $Y$ along
$Y'\times_Y Y''$ separates strict transforms of $Y'$ and $Y''$, as stated.
\end{proof}

\subsection{Regularization of strict transform}
\label{regstrsec} In this section we assume that $X$ is an integral noetherian scheme of dimension $d$ and $Y$
is a reduced closed subscheme whose maximal points are regular points of $X$ of codimension $1$. Note that
$(X,Y)_\sing$ does not contain maximal points of $Y$. Our aim is to prove the following statement.

\begin{prop}
\label{regprop} Assume that there is semi-strict embedded resolution of singularities over $X$ up to dimension
$<d$. Then there exists a blow up $f:X'\to X$ supported on $T=(X,Y)_\ssing$ such that the strict transform of
$Y$ is disjoint from $(X',f^{-1}(Y))_\ssing$.
\end{prop}

\begin{rem}\label{basarem}
It seems that the statement of proposition \ref{regprop} should hold true for strict desingularizations and
singular loci $(X,Y)_\sing$. Having such a result would allow to replace semi-strictness with strictness in
theorems \ref{mainth} and \ref{formmainth}, and proposition \ref{standardprop}.
\end{rem}

We need to track the behavior of both strict and total transforms of $Y$ (recall that the latter is the entire
preimage of $Y$) with respect to blow ups, so it will be more convenient to consider a more general situation.
In the sequel, $Z$ will be a scheme which remembers the history of total transforms and $T$ will denote a closed
set which we are allowed to modify. So, let $X$ and $Y$ be as above and $Z$ be a closed subscheme of $X$ which
contains $Y\cap(X,Y)_\ssing$ and is disjoint from $Y^0$. Note that $Y\cap Z$ is nowhere dense in $Y$,
$Y\setminus Z$ is a strictly monomial divisor in $X\setminus Z$, and $Y\cap(X,S)_\ssing\subset Z$ for the closed
set $S=Y\cup Z$. Let $T$ be a closed set with $Y\cap (X,S)_\ssing\subseteq T\subseteq Y\cap Z$. For any
$T$-supported blow up $f:X'\to X$ we use the following notation: $Y'$ is the strict transform of $Y$ in $X'$
(note that the morphism $Y'\to Y$ is birational), $Z'=Z\times_X X'$ and $S'=f^{-1}(S)=Y'\cup Z'$. Then the
proposition follows from the following more general lemma (take $Z=T$ in the proposition).

\begin{lem}
\label{reglem} Keep the above notation and assume that there is semi-strict embedded resolution of singularities
over $X$ up to dimension $<d$. Then there exists a $T$-supported blow up $f:X'\to X$ such that $Y'\subset
(X',S')_\sreg$.
\end{lem}
\begin{proof}
A required blow up will be obtained as a composition of few blow ups, which will gradually improve the strict
transform of $Y$. Note that while proving the lemma, we can replace $X$ with a neighborhood $X_0$ of $Y$ and
shrink $Z$ accordingly (i.e. replace $Z$ with $Z\cap X_0$). Indeed, if a $T$-supported blow up
$f_0:\Bl_R(X_0)\to X_0$ satisfies the assertion of the lemma for $X_0,Y$ and $Z_0$, then $f:\Bl_R(X)\to X$ is a
blow up of $X$ which extends $f_0$ trivially, and hence satisfies the assertion of the lemma (we use here that
$R$ is closed in $T\into X_0$, hence $R$ is closed in $X$).

Step 0. {\sl Given a $T$-supported blow up $f:X'\to X$, we can replace $X,Y,Z$ and $T$ with $X',Y',Z'$ and any
$T'$ with $Y'\cap (X',S')_\ssing\subseteq T'\subseteq Y'\cap f^{-1}(T)$.} First of all, we note that
$Y'\setminus Z'\toisom Y\setminus Z$ is a strictly monomial divisor in $X'\setminus Z'\toisom X\setminus Z$,
hence $X',Y',Z'$ and $T'$ satisfy the assumptions of the lemma. Suppose that the proposition holds for
$X',Y',Z'$ and $T'$, and let $f':X''\to X'$ be a $T'$-supported blow up with $Y''\subset (X'',S'')_\sreg$, where
$Y''$ is the strict transform of $Y'$ and $S''=f'^{-1}(S')$. The morphism $f''=f\circ f'$ is a composition of
$T$-supported blow ups, hence it is a $T$-supported blow up by lemma \ref{compblow}. Obviously, $Y''$ is the
strict transform of $Y$ in $X''$ and $S''$ is the preimage of $S$. Hence $f''$ solves our problem for $X,Y,Z$
and $T$.

Step 1. {\sl We can assume that $Y$ is irreducible.} Let $m$ be the number of irreducible components of $Y$. By
induction we can assume that $m>1$ and the lemma holds when $Y$ has less than $m$ irreducible components. Let
$Y=Y_1\cup Y_2$, where $Y_1$ is an irreducible component of $Y$ and $Y_2$ is the union of the others. The idea
is to construct a required blow up $X'\to X$ in two steps: achieve first that $Y_1\subset(X,S)_\sreg$, then
apply the induction assumption to $Y_2$. Let us check the details. Find a blow up $f:X'\to X$ which solves our
problem for $Y_1$, $Z_1=Y_2\cup Z$ and $T_1=Y_1\cap(X,S)_\ssing$. Obviously, $f$ is $T$-supported, so it
suffices to solve our problem for $X'$, $Y'$, $Z'$ and $T'=Y'\cap(X',S')_\ssing$.

Let $Y'_i$ denote the strict transforms of $Y_i$, then $S'$ is a strictly monomial divisor in a neighborhood of
$Y'_1$ and $S'=Y'_1\cup Y'_2\cup Z'$, in particular, $T'=Y'_2\cap(X',S')_\ssing$ is disjoint from $Y'_1$. Since
$Y'\setminus Z'$ is a strictly monomial divisor in $X'\setminus Z'$, we obtain that $Y'_2\setminus Z'$ is a
strictly monomial divisor as well. Thus, $X',Y'_2,Y'_1\cup Z'$ and $T'$ satisfy the assumptions the lemma. Since
$Y'_2$ contains $m-1$ irreducible components, there is a $T'$-supported blow up $f':X''\to X'$, which solves our
problem for $X',Y'_2,Y'_1\cup Z'$ and $T'$. Let us check that $f'$ solves our problem for $X',Y',Z'$ and $T'$
too. Indeed, $f'$ does not modify $Y'_1$ because $T'$ is disjoint from $Y'_1$. So, $S''=f'^{-1}(S')$ is a
strictly monomial divisors in neighborhoods of both $Y''_1\toisom Y'_1$ and $Y''_2$, and we obtain that it is a
strictly monomial divisor in a neighborhood of $Y''=Y''_1\cup Y''_2$.

We finished the only stage where a playing with $T$ is required. In the sequel, we automatically set $T'=Y'\cap
f^{-1}(T)$ for any $T$-supported blow up $f:X'\to X$ (i.e. $T'$ is chosen as large as possible).

Step 2. {\sl We can assume in addition to Step 1 that there exists a Cartier divisor $Y_2$ and a closed
subscheme $Y_1\into X$ such that $Y_2\setminus T=(Y\setminus T)\sqcup(Y_1\setminus T)$.} Note that $Y\cap
X_\reg$ is a Cartier divisor in $X_\reg$, hence $\Bl_Y(X)\to X$ is an isomorphism over $X\setminus(Y\cap
X_\sing)\supset X\setminus T$. Using lemma \ref{admisslem} we can find a $T$-supported blow up $f:X'\to X$
dominating $\Bl_Y(X)$, then $Y'_2=Y\times_X X'$ is a Cartier divisor. Define $Y'_1$ to be the schematic closure
of $Y'_2\setminus Y'$. Since $Y'\setminus f^{-1}(T)\toisom Y'_2\setminus f^{-1}(T)$ and $T'=Y'\cap f^{-1}(T)$,
we obtain that $Y'_2\setminus T'=(Y'\setminus T')\sqcup(Y'_1\setminus T')$. Replacing $X,Y,Z$ and $T$ with
$X',Y',Z'$ and $T'$, we achieve the condition of the step.

Step 3. {\sl We can strengthen the condition of Step 2 by achieving that $Y$ itself is a Cartier divisor.} Let
$\calI$ and $\calJ$ be the $\calO_X$-ideals of $Y_2$ and $Y\times_X Y_1$, respectively. By proposition
\ref{princprop} (iii), choosing sufficiently large $n$ and blowing up the ideal $\calI+\calJ^n$, we obtain a
$T$-supported blow up $f:X'\to X$ such that the strict transform of $Y$ is a Cartier divisor. Replace $X,Y,Z$
and $T$ with $X',Y',Z'$ and $T'$, as earlier. In the sequel, $\calI$ is the invertible ideal defining $Y$.

Step 4. {\sl We can assume in addition to Step 3 that $Y$ is regular, and $T$ and $W=Z\times_X Y$ are strictly
monomial divisors in $Y$.} For any point $x\in Y\setminus T$, there exists a neighborhood $U_x$ such that the
intersection of $S=Y\cup Z$ with $U_x$ is a strictly monomial divisor. Then $W\times_X U_x$ is a strictly
monomial divisor in $Y\times_X U_x$, and we therefore obtain that $(Y,W)_\ssing\subset T$. Since $\dim(Y)\le
d-1$ and, by the assumption of the lemma, there is embedded resolution of singularities over $X$ in dimensions
smaller than $d$, there exists a closed $T$-supported subscheme $R\into Y$ such that $\calY'=\Bl_R(Y)$ is
regular, $W'=W\times_Y\calY'$ is a monomial divisor in $\calY'$ and $\calT'=T\times_Y\calY'$ is a Cartier
divisor. Furthermore, by the following lemma blowing up self-intersections of $W'$ (which lie above $T$) we can
achieve that $W'$ is strictly monomial. Note that this operation does not destroy the other properties we have
established.

\begin{lem}\label{strnorm}
Given a regular scheme $X=X_0$ with a normal crossing divisor $Z$, there exists a sequence of blow ups
$X_n\stackrel{f_n}\to\dots\stackrel{f_2}\to X_1\stackrel{f_1}\to X_0$ such that each $X_i$ is regular, each
$Z_i=f_i^{-1}(Z_{i-1})$ is normal crossing, $Z_n$ is strictly normal crossing, and the center of each $f_i$ is a
regular subscheme which is a self-intersection of $Z_{i-1}$ of maximal multiplicity. In particular, the
composite blow up $X_n\to X$ is supported over $(X,Z)_\ssing$.
\end{lem}
\begin{proof}
Let $\oZ_1\.\oZ_l$ be the irreducible components of $Z$, and choose any $\oZ=\oZ_i$ with non-empty
self-intersection. Let $T$ be the self-intersection of $\oZ$ of maximal multiplicity $n(\oZ)$ (i.e. each point
of $T$ has $n(\oZ)$ preimages in the normalization of $\oZ$). Then similarly to \cite[7.2]{dJ} one checks that
$T$ is a regular closed subscheme of $\oZ$ which is transversal to all other components of $Z$. Blowing up $X$
along $T$ we obtain a regular scheme $X_1$ such that $Z_1$ is normal crossing and the preimage of $\oZ$ consists
of two components: a regular exceptional component, and the strict transform $\oZ'$ of $\oZ$. Since $\oZ'$ is
isomorphic to the blow up of $Z$ along $T$, we obtain that $n(\oZ')<n(\oZ)$. Now it is clear, that we can
iterate the same process by picking up any irreducible component of $Z_1$ with non-empty self-intersection, and
the process will stop after $n\le \sum_{i=1}^l n(\oZ_i)$ steps. Then $X_n\to X$ is as required, and clearly we
only modified $X_i$'s over the set $(X,Z)_\ssing$ where the normal crossing divisor $Z$ is not strict.
\end{proof}

Consider $R$ as a closed subscheme of $X$, and let $\calJ$ be its $\calO_X$-ideal. By proposition
\ref{princprop} (ii), there exists $n$ such that the strict transform of $Y$ in $X'=\Bl_{\calI+\calJ^n}(X)$ is a
Cartier divisor. Define $Y',Z'$ and $T'$ as usual, then $Y'\toisom\calY'$ by \ref{princprop} (i). In particular,
$\calT'\toisom (T\times_X X')\times_{X'}Y'$ and we obtain that $|\calT'|=T'$. To check that $X',Y',Z'$ and $T'$
satisfy the conditions of the step, we note that $Z'\times_{X'}Y'\toisom Z\times_X Y'\toisom W\times_X Y'\toisom
W'$ is a strictly monomial divisor in $Y'$. Finally, $T'$ is a divisor supported on $W'$, hence it is a strictly
monomial divisor too.

Step 5. {\sl We can achieve in addition to Step 4 that $X$ is regular.} Note that $X$ is regular in a
neighborhood of $Y$ because $Y$ is a regular Cartier divisor. So, we can simply shrink $X$.

Step 6. {\sl We can achieve in addition to Step 5 that $Z$ is a Cartier divisor.} Let $D$ be the divisorial part
of $Z$, i.e. the schematic closure of $\sqcup_{z\in Z\cap X^1}\Spec(\calO_{Z,z})$ in $X$. Also, let
$\calI_Z\subset\calI_D\subset\calO_X$ be the ideals of $D$ and $Z$. Since $X$ is regular, $\calI_D$ is
invertible and we obtain a splitting $\calI_Z=\calI_D\calI_\tilZ$ where $\calI_\tilZ$ is an ideal supported in
codimension at least two. The support of the scheme $\tilZ=\Spec(\calO_X/\calI_\tilZ)$ is the locus of $Z$ where
it is not a divisor (large codimension or embedded components), hence $\tilZ\cap Y\subset T$. Now, blowing up
$X$ along $\tilZ\times_X Y$ we obtain a $T$-supported blow up $X'\to X$ such that the ideal
$\calI_\tilZ\calO_{X'}$ is principal in a neighborhood of $Y'$ by lemma \ref{strcarlem}. Then it is clear that
the closed subscheme given by the ideal $\calI_D\calI_\tilZ\calO_{X'}$ is principal in that neighborhood of $Y'$
as well. Thus, we can achieve that $Z$ is a Cartier divisor at cost of possible destroying the conditions of
Steps 2--5. Since the property of $Z$ being a Cartier divisor is preserved by any modification of $X$, we should
simply rerun Steps 2--5 once again.

The remaining part of the proof is more or less standard: it will suffice only to blow up some components of $T$
(which are regular subschemes of codimension $2$ in a regular scheme $X$). We prefer to give a detailed proof
mainly for the sake of completeness.

Step 7. {\sl Let $T_1\. T_n$ be the irreducible components of $T$, then we can achieve in addition to Step 6
that each $T_i$ belongs to a unique irreducible component $Z_i$ of $Z$, $T_i=Z_i\cap Y$ and $Z_i\cup Y$ is a
strictly monomial divisor.} Consider $T_1$ as a reduced closed subscheme, and let $\calJ\subset\calO_X$ be its
ideal and $m$ be its multiplicity in $W$. Set $X'=\Bl_{\calI+\calJ^m}(X)$ and define $Y'$,$Z',T'$ as usual, then
$Y'\toisom \Bl_{T_1}(Y)\toisom Y$. Note that $W'=Y'\times_{X'}Z'\toisom Y'\times_X Z\toisom Y'\times_Y W$ is
isomorphic to $W$, hence the conditions of Step 4 are still satisfied. It follows from the next lemma that only
one component of $Z'$, say $Z_1$, contains $T_1$, $T_1=Z_1\cap Y'$ and $Z_1\cup Y'$ is strictly monomial in a
neighborhood of $Y'$. Hence shrinking $X'$ (i.e. replacing it with $X_1$ in the notation of the lemma) and
replacing $X,Y,Z,T$ with $X',Y',Z',T'$, we achieve that $T_1$ satisfies the conditions of Steps 5--7. It remains
then to repeat this procedure for other $T_i$'s and to note that blowing up $T_2$ preserves monomiality of
$Z_1\cup Y$ by part (iii) of the lemma, similarly for $T_3$, etc.

\begin{lem}
Assume that $X$ is a regular scheme, $Y=\Spec(\calO_X/\calI)\into X$ is a regular divisor,
$T=\Spec(\calO_X/\calJ)\into Y$ is a regular divisor in $Y$ and $m$ is a positive natural number. Consider the
blow up $f:X'=\Bl_{\calI+\calJ^m}\to X$ and let $Y'$ be the strict transform of $Y$, then

(i) in a neighborhood of $Y'$, $Y\times_X X'$ is a strictly monomial divisor which factors as $Y'\cup m\tilT$,
where $\tilT$ is the preimage of $T$ with the induced reduced scheme structure;

(ii) if $Z$ is a Cartier divisor in $X$ with $Y\times_X Z=mT$, then $Z\times_X X'$ coincides with $Y'\cup
m\tilT$ in a neighborhood of $Y'$;

(iii) if $\tilZ$ is a divisor in $X$ such that $Y\cup\tilZ$ is strictly monomial and $(Y\times_X\tilZ)\cup T$ is
strictly monomial in $Y$, then $Y'\cup (\tilZ\times_X X')$ is strictly monomial in a neighborhood of $Y'$.
\end{lem}
\begin{proof}
The statement is local on $X$, so we can assume that $X=\Spec(A)$, $I=\calI(X)=(x)$ and $J=\calJ=(x,y)$. Then
$\calI+\calJ^m$ corresponds to the ideal $L=(x,y^m)$ and $X'$ is pasted from the charts
$X_1=\Spec(A[\frac{L}{y^m}])$ and $X_2=\Spec(A[\frac{L}{x}])$. The strict transform of $Y$ is disjoint from
$X_2$, hence we can restrict our study to $X_1$, and we will actually show that it is a required neighborhood of
$Y'$. Note that the $A$-algebra $B=A[S]/(y^mS-x)$ defines a regular subscheme in $\bfA^1_X$, in particular, $B$
has no $y^m$-torsion and therefore  the surjection $B\to A[\frac{L}{y^m}]$ is an isomorphism.  It follows that
$X_1\toisom\Spec(B)$ is regular and $Y\times_X X_1$ is isomorphic to the strictly monomial subscheme of
$\Spec(B)$ given by the condition $y^mS=0$, as stated in (i).

Since the divisor $Y\times_X X_1$ in $X_1$ is defined by $y^mS=0$, it coincides with $Y'\cup m\tilT$, where
$\tilT$ is given by $y=0$ (i.e. $\tilT$ is the set-theoretical preimage of $T$ in $X_1$). In particular,
$T':=Y'\times_{X_1}\tilT$ is the zero locus of $(S,y)$, hence it coincides with the preimage of $T$ in $Y'$.
Let, now, $Z$ be as in (ii) and $Z_1=Z\times_X X_1$. Note that $Y'\toisom\Bl_{\calJ^m}(Y')\toisom Y'$ and hence
$T'\toisom T$. Therefore, $Y'\times_{X_1}Z_1\toisom Y\times_X Z=mT\toisom mT'$. Since $m\tilT$ is an irreducible
component of $Z_1$ and its intersection with $Y'$ is $mT'$, we obtain that $Y'$ and $m\tilT$ are the only
irreducible components of $Z_1$ that are not disjoint from $Y'$, so we obtain (ii). Finally, (iii) is proved by
an explicit computation similar to the proof of (i), so we omit the details.
\end{proof}

Step 8. {\sl We can achieve in addition to Step 7 that for any irreducible component $\tilZ$ of $Z$ the divisor
$\tilZ\cup Y$ is strictly monomial.} Let us prove first that any two irreducible components $\tilW,\tilW'$ of
$\tilZ\cap Y$ are disjoint. Assume, on the contrary, that $V=\tilW\cap\tilW'$ is non-empty. Since $W$ is a
strictly monomial divisor, we obtain that $V$ is of codimension $2$ in $Y$ and $\tilW,\tilW'$ are the only
components of $W$ which contain $V$. Note that $\tilW$ is not contained in $T$ because otherwise Step 7 would
imply that $\tilW=\tilZ\cap Y$. By the same reason, $\tilW'$ is not in $T$ and, since $T$ is a union of
components of $W$ by Step 4, we obtain that $V\setminus T$ is not empty, say contains a point $y$. But the
latter is an absurd because $Z\cup Y$ is strictly monomial at $y$, but $y$ belongs to two different irreducible
components of $\tilZ\cap Y$. The contradiction shows that $V$ is actually empty.

Now, we are ready to check that $\tilZ\cup Y$ is a strictly monomial divisor in a neighborhood of $Y$. Since $Y$
itself is strictly monomial it is enough to check that $\tilZ\cup Y$ is a strictly monomial divisor in a
neighborhood of any point $y\in\tilZ\cap Y$. Let $\tilW$ be the unique irreducible component of $\tilZ\cap Y$
that contains $y$. We can assume that $\tilW$ is not contained in $T$ as the latter case was dealt with in Step
7. Shrinking $X$ (and $Y$) we can assume that $\tilZ\times_X Y=\tilW$. The irreducible divisor $\tilZ$ is of the
form $m\tilZ'$ where $\tilZ'$ is reduced. Note that $\tilZ\cup Y$ is a strictly monomial divisor in a
neighborhood of the generic point $\eta\in\tilW$ because $\tilW$ is not contained in $T$. It follows that
$\tilW'=\tilZ'\times_X Y$, which is an irreducible divisor in $Y$, is reduced at $\eta$. Therefore, $\tilW'$ is
an integral divisor in $Y$, and actually $\tilW'=\tilW_\red=\frac{1}{m}\tilW$.

It suffices to show that $\tilZ'\cup Y$ is strictly monomial, so we can replace $\tilZ$ with its reduction
achieving that $\tilW$ becomes reduced. We can check monomiality locally at each point $x\in\tilW$. There exists
a regular sequence of parameters $x_1\. x_n\in\calO_{X,x}$ such that $Y=V(x_1)$ and $\tilW=V(x_1,x_2)$ locally
at $x$ (we use that $X$ is regular, $Y$ is a regular divisor in $X$ and $W$ is a regular divisor in $Y$). Since
$\tilZ$ is a divisor, it is of the form $V(f)$, and then the image $f'\in\calO_{X,x}/(x_1)$ of $f$ is of the
form $u'x_2$ with a unit $u'\in\calO_{X,x}/(x_1)$ because locally at $x$ $f'$ defines the closed subscheme
$\tilZ\times_X Y=W$ in $Y$. Lifting $u'$ to a unit $u\in\calO_{X,x}$ and replacing $f$ with $f/u$ we can get rid
of $u$ and $u'$, and then $f=x_2+x_1y$ for some $y\in\calO_{X,x}$. In particular, it becomes clear that
$\tilZ\cup Y=V(x_1(x_2+yx_1))$ is a strictly monomial divisor locally at $x$. By compactness of $Y$, $\tilZ\cup
Y$ is strictly monomial in a neighborhood of $Y$, hence it suffices simply to shrink $X$.

Step 9. {\sl If the conditions of Step 8 are satisfied, then $S$ is a strictly monomial divisor in a
neighborhood of $Y$.} Let $Z_1\. Z_n$ be the irreducible components of $Z$ which have non-empty intersection
with $Y$. By the previous step, $Y\cup Z_i$ is a strictly monomial divisor for any $i$. Since $\cup_i(Y\cap
Z_i)$ is a strictly monomial divisor in $Y$, it follows that $Y\cup (\cup_i Z_i)$ is a strictly monomial divisor
in a neighborhood of $Y$ (replace $Z$ with its reduction, then for any point $y\in Y$ the claim reduces to
linear algebra in the tangent space $T_y$). It finishes the proof of Step 9 and concludes the proof of the
theorem.
\end{proof}

\subsection{Main results}\label{mainsec}

\begin{lem}\label{lemlem}
Let $X$ be an integral quasi-excellent scheme of characteristic zero with a closed subset $Z$ such that
$T=(X,Z)_\ssing$ is of finite type over a field $k$. Assume that $\dim(X)=d$ and there is semi-strict embedded
resolution of singularities over $X$ up to dimension $<d$. Then the pair $(X,Z)$ admits a semi-strict
desingularization.
\end{lem}
\begin{proof}
Note that for any $T$-supported blow up $f:X'\to X$ with $Z'=Z\times_X X'$, the scheme $T'=(X',Z')_\ssing$ is of
finite type over $k$ because it sits over $T$, and hence is of finite type over $T$. In particular, while
proving the lemma we will freely replace $X,Z $ and $T$ with $X',Z'$ and $T'$ as above. Let $f:X'\to X$ be a
$T$-supported blow up dominating $\Bl_Z(X)\to X$; it exists by lemma \ref{admisslem}. Then $Z'=f^{-1}(Z)$ is the
support of the Cartier divisor $Z\times_X X'$. Replace $X$ and $Z$ with $X'$ and $Z'$. By corollary
\ref{desingcor} applied to $X$ and $T\supset X_\sing$, there exists a strict desingularization $f:X'\to X$ of
$X$. Replacing $X$ and $Z$ with $X'$ and $f^{-1}(Z)$, we achieve in addition that $X$ is regular. Now, $X$ and
$Y=Z$ satisfy assumptions of proposition \ref{regprop}, hence there exists a $T$-supported blow up $f:X'\to X$
such that the strict transform $Y'$ of $Z$ is disjoint from $T'=(X',S')_\ssing$ for $S'=f^{-1}(Z)$.

The Zariski closure $Z'$ of $S'\setminus Y'$ is of finite type over $k$ because $Z'\subset f^{-1}(T)$. Note that
$(X',Z')_\ssing\subset T'\subset Z'$, hence we can apply corollary \ref{desingcor} to find a strict
desingularization $g:X''\to X'$ of the pair $(X',Z')$. Let $Y'',Z''$ and $S''$ be the preimages of $Y',Z'$ and
$S'$. Since $g$ is $T'$-supported, $T'\cap Y'=\emptyset$ and $S'$ is a monomial divisor in a neighborhood of
$Y'$, we obtain that $S''$ is a monomial divisor in a neighborhood of $Y''\toisom Y'$. Also, $S''=Y''\cup Z''$
and $Z''$ is a monomial divisor, hence the $S''$ is a monomial divisor. The induced morphism $f':X''\to X$ is a
$T$-supported blow up because $T'\subset Z'\subset f^{-1}(T)$. Therefore $f'$ provides a semi-strict
desingularization of $(X,Z)$.
\end{proof}

Now, we are prepared to prove theorem \ref{mainth}.
\begin{proof}[Proof of theorem 1.1]
As was mentioned in the introduction, it suffices to prove that there is semi-strict embedded resolution of
singularities over a quasi-excellent noetherian scheme $k$ of characteristic zero. We will prove by induction on
$d$ that there is semi-strict embedded resolution of singularities over $k$ up to dimension $<d$. The case $d=0$
is trivial because $X^{<0}=\emptyset$. Assume that there is semi-strict embedded resolution of singularities
over $k$ up to dimension $<d-1$. By proposition \ref{locdesprop}, it suffices to prove that, if $S$ is a local
integral $k$-scheme of essentially finite type over $k$ and dimension $d$, $s\in S$ is the closed point,
$f:S'\to S$ is a modification and $Z'\subset S'$ is a closed subset such that $T'=(S',Z')_\sing$ is contained in
$f^{-1}(s)$, then $(S',Z')$ admits a semi-strict desingularization.

Set $R=(S',Z')_\ssing$. We claim that there exists an $R$-supported blow up $g:S''\to S'$ with $Z''=g^{-1}(Z')$
such that $(S'',Z'')_\ssing$ is contained in the preimage of $s$. Indeed, by lemma \ref{strnorm}, there exists a
blow up $\tilg:\tilS''\to \tilS':=(S',Z')_\reg$ which is supported on $(\tilS',\tilZ')_\ssing$, where
$\tilZ'=Z'\cap\tilS'$, and such that $\tilS''$ is regular and the preimage $\tilg^{-1}(\tilZ')$ of the monomial
divisor $\tilZ'$ is strictly monomial (we use that the latter happens iff the preimage of the closed set
$|\tilZ'|$ is a strictly normal crossing divisor). Extending $\tilg$ to a blow up $g:S''\to S'$ we obtain an
$R$-supported blow up $g:S''\to S'$ with $Z''=g^{-1}(Z')$ such that $(S'',Z'')_\sreg$ contains the preimage of
$(S',Z')_\reg$. It follows that $(S'',Z'')_\ssing$ is contained in the preimage of $s$, which is a
$k(s)$-variety.

It suffices to find a semi-strict desingularization of the pair $(S'',Z'')$ because any such desingularization
is also a semi-strict desingularization of the original pair $(S',Z')$. By proposition \ref{locdesprop}, there
is semi-strict embedded resolution of singularities over $S$ up to dimension $<d$ because any local $S$-scheme
of essentially finite type is also a $k$-scheme of essentially finite type. In particular, there is semi-strict
embedded resolution of singularities over $S''$ up to dimension $<d$, hence the pair $(S'',Z'')$ admits a
semi-strict desingularization by lemma \ref{lemlem}.
\end{proof}

We will deduce desingularization of formal schemes. By a desingularization of a formal pair $(\gtX,\gtZ)$ we
mean a formal blow $\gtX'\to\gtX$ supported on $\gtX_\sing\cup\gtZ_\sing$ and such that
$\gtX'=(\gtX',\gtZ')_\reg$, where $\gtZ'=\gtZ\times_\gtX\gtX'$. We have already seen in \S\ref{formregsec} that
it is not so easy to define quasi-excellent formal schemes, and now we have to introduce one more notion. We say
that $\gtX$ is a {\em universally quasi-excellent} formal scheme if any formal $\gtX$-scheme of finite type is
quasi-excellent. Obviously, any special formal scheme is universally quasi-excellent (and, as we noted in remark
\ref{Gabrem}, Gabber proved that any quasi-excellent formal scheme is universally so). Consider such a formal
scheme $\gtX$ with a closed subscheme $\gtZ$.

\begin{cor}
\label{formcor} Assume that $\gtX=\hatBl_J(\gtX_0)$ and $\gtZ=\gtZ_0\times_{\gtX_0}\gtX$, where $\gtX_0=\Spf(A)$
is a reduced universally quasi-excellent formal scheme of characteristic zero, $I,J\subset A$ are ideals and
$\gtZ_0=\Spf(A/I)$. Then the pair $(\gtX,\gtZ)$ admits a desingularization.
\end{cor}
\begin{proof}
Set $X_0=\Spec(A)$, $X=\Bl_J(X_0)$ and $Z=\Spec(A/I)\times_{X_0}X$, then the pair $(\gtX,\gtZ)$ is isomorphic to
the $P$-adic completion of the pair $(X,Z)$, where $P$ is an ideal of definition of $A$. By the previous
theorem, the pair $(X,Z)$ admits a desingularization $X'\to X$ (use lemma \ref{dummylem} if $X$ is not
integral). Using lemma \ref{compllem}, one checks that the $P$-adic completion $\gtX'\to\gtX$ of the blow up
$X'\to X$ is a required desingularization of the pair $(\gtX,\gtZ)$.
\end{proof}

One could expect that the corollary allows to desingularize an arbitrary universally quasi-excellent formal
scheme of characteristic zero by patching local desingularizations (proposition \ref{locdesprop} does such
patching job in the case of schemes). Unfortunately, it cannot be done in general because not open ideals does
not have to extend from an open formal subscheme. For this reason we are forced to consider the case when
blowing up an open ideal suffices for desingularization, i.e. the case then $\gtX$ is rig-regular and $\gtZ$ is
rig-monomial.

Recall that a (formal) scheme $X$ is called {\em quasi-paracompact} if it admits a covering $\{X_i\}_{i\in I}$
of {\em finite type} (i.e. each $X_i$ intersects only finitely many $X_j$'s) with open quasi-compact $X_i$'s.
Any irreducible quasi-paracompact locally noetherian scheme is actually noetherian, but quasi-paracompactness is
a much more interesting property in the case of formal schemes. For example, Drinfeld's upper half plane,
non-Archimedean Stein spaces and analytifications of varieties over a non-Archimedean field admit
quasi-paracompact formal models, which can be chosen to be irreducible if the corresponding non-Archimedean
space is irreducible. (Irreducibility is understood here in the sense of \cite{Con2}; it is a rather subtle
notion because it is not preserved by localizations, unlike the scheme case.)

Note that Berkovich considered in \cite{Ber2} (and later works) quasi-para\-compact formal schemes of locally
finite presentation over the ring of integers of a non-Archimedean field (they are simply called formal schemes
of locally finite presentation in loc.cit.). One can define analytic generic fiber for such formal schemes.
Recently, Bosch extended Raynaud's theory to quasi-paracompact formal schemes and rigid spaces, see
\cite[2.8.3]{Bo}. Let us say that a (formal) scheme is para-noetherian if it is quasi-para\-compact and locally
noetherian.

\begin{theor}
\label{formmainth} Let $\gtX$ be a reduced universally quasi-excellent para-Noethe\-rian formal scheme of
characteristic zero and $\gtZ$ be a closed formal subscheme. Assume that $\gtX$ is rig-regular and $\gtZ$ is a
rig-regular divisor, then the pair $(\gtX,\gtZ)$ admits a desingularization.
\end{theor}
\begin{proof}
The proof exploits the same reasoning as was used in the proofs of proposition \ref{locdesprop} and theorem
\ref{formtheor}. Since there are mild complications due to lack of quasi-compactness, we give the full argument.
By our rig-assumptions, $\gtT=\gtX_\sing\cup\gtZ_\sing$ is a reduced closed subscheme of $\gtX_s$. It follows
that $(\gtX',\gtZ\times_\gtX\gtX')_\sing$ is an $\gtX_s$-supported scheme for any formal blow up
$\gtX'=\hatBl_\gtI(\gtX)$ with open $\gtI$.

Each connected component can be desingularized separately, so assume that $\gtX$ is connected. Choose a locally
finite open affine covering $\{\gtX_i\}_{i\in I}$. By connectedness of $\gtX$, $I$ is at most countable. We
assume that $I=\bfN$ because the case of finite $I$ is similar and easier. By $\gtU_i$ we denote the open set
$(\gtX\setminus\gtT)\cup(\cup_{j\le i}\gtX_j)$. Let us say that an ideal $\gtI\subset\calO_\gtX$ has a
quasi-compact support if it is trivial outside of a quasi-compact formal open subscheme. We will use induction
whose step is the following statement. {\sl Let $\gtI$ be a $\gtT$-supported open ideal with a quasi-compact
support such that the blow up $f:\gtX'=\hatBl_\gtI(\gtX)\to\gtX$ desingularizes the pair $(\gtX,\gtZ)$ over
$\gtU_{i-1}$. Then there exists an open $\gtT$-supported ideal $\gtL\subset\gtI$ with a quasi-compact support
such that $\gtL|_{\gtU_{i-1}}=\gtI|_{\gtU_{i-1}}$ and $\hatBl_\gtL(\gtX)\to\gtX$ desingularizes $(\gtX,\gtZ)$
over $\gtU_i$.}

Assume the induction step for now. Then by induction we can find a sequence of ideals
$\calO_\gtX=\gtI_0\supset\gtI_1\supset\dots$ such that the sequence stabilizes on each $\gtX_i$ starting with
$\gtI_i$, and each morphism $\Bl_{\gtI_i}(\gtX)\to\gtX$ desingularizes $(\gtX,\gtZ)$ over $\gtU_i$. Now, it is
clear that one can desingularize $(\gtX,\gtZ)$ by blowing up the ideal $\gtI=\cap\gtI_i$.

It remains to establish the induction step. Recall that $\gtX_i$ is affine and set
$\gtU'_{i-1}=\gtU_{i-1}\times_\gtX\gtX'$, $\gtX'_i=\gtX_i\times_\gtX\gtX'$ and $\gtZ'_i=\gtZ\times_\gtX\gtX'_i$.
The pair $(\gtX'_i,\gtZ'_i)$ admits a desingularization $\hatBl_{\gtJ_0}(\gtX'_i)\to\gtX'_i$ by corollary
\ref{formcor}. Note that $(\gtX'_i)_\sing\cup(\gtZ'_i)_\sing$ is a closed subscheme of $(\gtX'_i)_s$ disjoint
from $\gtU'_{i-1}$, hence the ideal $\gtJ_0$ is $\gtT'_i$-supported for the closed subset
$\gtT'_i=\gtX'_i\setminus\gtU'_{i-1}$ of $\gtX'_i$, in particular, $\gtJ_0$ is an open ideal. Choose a
quasi-compact neighborhood $\ogtX'$ of the Zariski closure $\ogtT'_i$ of $\gtT'_i$. By lemma
\ref{formextendlem}, $\gtJ_0$ extends to a $\ogtT'_i$-supported ideal $\ogtJ\subset\calO_{\ogtX'}$. Since the
support of $\ogtJ$ is closed in $\gtX'$, it can be extended trivially to a $\ogtT'_i$-supported ideal
$\gtJ\subset\calO_{\gtX'}$. In particular, $\gtJ$ is an open $\gtT$-supported ideal with a quasi-compact
support.

By lemma \ref{formsupplem}, the morphism $f':\gtX''=\Bl_\gtJ(\gtX')\to\gtX$ is isomorphic to a $\gtT$-supported
blow up $\hatBl_\gtL(\gtX)$ with a quasi-compact support (the lemma is formulated for noetherian formal schemes,
but the same proof works for para-noetherian formal schemes and $\gtT$-supported blow ups with a quasi-compact
support). Replacing $\gtL$ with $\gtL\gtI$ we achieve that $\gtL$ is contained in $\gtI$ (note that
$\hatBl_{\gtL\gtI}(\gtX)\toisom\gtX''$ because $\gtX''$ dominates $\gtX'=\hatBl_\gtI(\gtX)$). It remains to
check that $\gtL$ is as required. The ideal $\gtL$ is open because it is $\gtT$-supported. Next, $f'$
desingularizes $(\gtX,\gtZ)$ over $\gtU_{i-1}$ because it is isomorphic to $f$ over it (the blow up
$\gtX''\to\gtX'$ is $\ogtT'_i$-supported, hence it is an isomorphism over $\gtU'_{i-1}$). Finally, $f'$
desingularizes $(\gtX,\gtZ)$ over $\gtX_i$ because $\gtX''\times_\gtX\gtX_i\toisom\hatBl_{\gtJ_0}(\gtX'_i)$.
\end{proof}

It seems natural to expect that the rig-assumptions in the above theorem can be removed, but perhaps one has to
use stronger methods to prove such a generalization. Some kind of canonical desingularization should be used
because a badly chosen desingularization of an open formal subscheme can have no extension to the entire formal
scheme.

\appendix

\section{Standard desingularization}
The aim of this appendix is to prove that embedded desingularization in our sense can be deduced from a
desingularization result of Hironaka's type.

\begin{definsect}
\label{standardrem} We say that there is {\em standard resolution of singularities up to dimension $<d$} over a
locally noetherian scheme $k$ if there is resolution of singularities up to dimension $<d$ over $k$, and for any
regular scheme $X$ of finite type over $k$, if $\dim(X)<d$, $E\subset X$ is a normal crossing divisor and
$Z\subset X$ is a closed subset, then there exists a $Z$-supported blow up $f:X'\to X$ such that $f^{-1}(E\cup
Z)$ is a normal crossing divisor.
\end{definsect}

The definition is motivated by an observation that the work of Hironaka and all recent desingularization works
imply standard desingularization. We refer the reader to \cite{Wl}, \S\S1, 2.1 and 2.3 for an excellent
exposition of the general strategy shared by all known desingularization proofs, and give here only a very brief
explanation. Recall that the {\em order} or {\em multiplicity} $\mu_x$ of an ideal $\calI\subset\calO_X$ (or the
corresponding closed subscheme) at $x\in X$ is the maximal number $n$ for which $\calI_x\subseteq m_x^n$. If $n$
is the maximal multiplicity of $\calI$ at a point, $f:X'\to X$ is a blow up of $X$ along a regular subvariety
lying in the multiplicity $n$ locus of $\calI$ and $E'$ is the exceptional divisor, then $\calI\calO_{X'}(nE')$
is an ideal in $\calO_{X'}$ of maximal multiplicity $n$, which is called the {\em weak} (or controlled)
transform of $\calI$ under $f$. (Actually, we take the full transform $\calI\calO_{X'}$ and factor out an
obvious divisorial part. So, the weak transform lies somewhere on the way from the full to the strict transform,
but unlike the strict transform it can be easily described.)

The main ingredient of desingularization proofs is the following statement: let $X,E,Z$ be as above, and assume
that the multiplicity of $Z$ at the points of $X$ is at most $\mu$, then there is a composition $g:X'=Y_n\to
Y_{n-1}\to\dots\to Y_0=X$ of blow ups with regular centers which are contained in the maximal multiplicity loci
of the weak transforms of $Z$ and such that the union of the exceptional divisor $E_g$ with $g^{-1}(E)$ is a
normal crossing divisor and the multiplicity of the weak transform of $Z$ at the points of $X'$ is at most
$\mu-1$. For example, the above statement is Main Theorem II in \cite{Hir}, or resolution of marked ideals in
\cite[2.1.3]{Wl}. Applying this procedure $\mu$ times, we obtain a composition of $Z$-supported blow ups
$f:X_0\to X_1\to\dots\to X_\mu=X$ such that the union of the exceptional divisor $E_f$ with $f^{-1}(E)$ is a
normal crossing divisor and the weak transform of $Z$ is empty. Then $E_f$ contains the preimage of $Z$ and,
therefore, $f^{-1}(E\cup Z)$ is a normal crossing divisor. Thus, standard desingularization of algebraic
varieties of characteristic zero follows from \cite{Hir}, \cite{Vil}, \cite{BM}, \cite{Wl} and other
desingularization works.

\begin{propsect}
\label{standardprop} Let $k$ be a noetherian scheme. If there is standard resolution of singularities over $k$
up to dimension $<d$, then there is semi-strict embedded resolution of singularities over $k$ up to dimension
$<d$.
\end{propsect}
\begin{proof}
Using induction on $d$ we can assume that there is semi-strict resolution of singularities over $k$ up to
dimension $<d-1$. By proposition \ref{locdesprop}, it suffices to show that if $X$ is a blow up of a local
$k$-scheme of dimension $d-1$, in particular, $\dim(X)\le d-1$, and $Z\into X$ is a closed subscheme, then the
pair $(X,Z)$ admits a semi-strict desingularization.

Step 1. {\sl We can assume that $X$ is regular and $Z$ is a reduced Cartier divisor.} Set $V=(X,Z)_\sreg$. As we
saw in the beginning of the proof of \ref{locdesprop}, the blow up $\Bl_Z(X)\to X$ is dominated by a
$V$-admissible blow up $X'\to X$, and replacing $(X,Z)$ with $(X',Z\times_X X')$, we achieve that $Z$ is a
Cartier divisor. Since we assume standard resolution of singularities over $k$ up to dimension $<d$, there
exists an $X_\reg$-admissible blow up $X'\to X$ with regular $X'$. Replacing $(X,Z)$ with $(X',Z\times_X X')$,
we can assume that $X$ is regular. In particular, it is now harmless to replace $Z$ with its reduction.

In the sequel, we will desingularize $(X,Z)$ by a sequence of blow ups. Let $T$ denote the set we are allowed to
modify. Clearly, we have to start with $T=(X,Z)_\ssing$, but it will not be so in the sequel.

Step 2. {\sl We can assume in addition to Step 1 that $Z=Y\cup T$, where $Y$ is a divisor disjoint from
$(X,Z)_\ssing$.} (We warn the reader that $Y=\emptyset$ does not have to work fine because $T$ can be strictly
smaller than $Z$.) By proposition \ref{regprop}, there exists a $T$-supported blow up $f:X'\to X$ with the
following property: if $Y'$ is the strict transform of $Z$ and $Z'=f^{-1}(Z)$, then $Y'\subset(X',Z')_\sreg$.
Note that $Z'=Y'\cup T'$, where $T'=f^{-1}(T)$, therefore $X',Y',Z'$ and $T'$ satisfy the claim of the step with
the only possible exception: it can happen that $X'$ is singular. Find a desingularization $f':X''\to X'$ and
set $Y''=f'^{-1}(Y')$, $Z''=f'^{-1}(Z')$ and $T''=f'^{-1}(T')$. Note that $f'$ is a $T$-supported blow up
because $X'_\sing\subset T'$. Also, $f'$ is an isomorphism over a neighborhood of $Y'$ because $Y'\subset
X'_\reg$. Hence $Y''\subset(X'',Z'')_\sreg$, and we can replace $X',Y',Z'$ and $T'$ with $X'',Y'',Z''$ and
$T''$, which are as claimed.

The rest is obvious. By the definition of standard desingularization, there exists a $T$-supported blow up
$f:X'\to X$ with monomial $f^{-1}(Z)$.
\end{proof}

\end{document}